\documentclass[11pt,a4paper]{article}

\usepackage{amsmath,amssymb}

\newcommand{\diam}{\operatorname{diam}}

\newcommand{\counte}{theorem}

\allowdisplaybreaks

\renewcommand{\thefootnote}{\fnsymbol{footnote}}

\begin{document}

\renewcommand{\thefootnote}{\arabic{footnote}}

\centerline{\bf Rigidity theorems for glued spaces being suspensions, cones and joins}
\centerline{\bf in Alexandrov geometry with curvature bounded below
\footnote{Supported by NSFC 11001015 and 11171025.
\hfill{$\,$}}}

\vskip5mm

\centerline{ Xiaole Su, Hongwei Sun, Yusheng Wang\footnote{The corresponding author (E-mail: wwyusheng@gmail.com). \hfill{$\,$}}}

\vskip6mm

\noindent{\bf Abstract.} In the paper, we give rigidity theorems when the glued space of two Alexandrov spapces with curvature bounded below
is a suspension, cone or join. And we list some
basic properties of joins in Appendix.

\vskip1mm

\noindent{\bf Key words.} Alexandrov spaces, gluing, suspensions, cones, joins

\vskip1mm

\noindent{\bf Mathematics Subject Classification (2000)}: 53-C20.

\vskip6mm

An important and interesting class of Alexandrov spaces with
curvature $\geq k$ is provided with nonempty boundary (cf. [BGP]). A
fundamental and significant result on such spaces is the Gluing
Theorem, which is formulated as follows ([Pet]):

\vskip1mm

{\it Let $M_i$ $(i=1,2)$ be Alexandrov spaces with curvature $\geq
k$ and nonempty boundary $\partial M_i$. Let there be an isometry
$i:\partial M_1\to\partial M_2$, where $\partial M_i$ are considered
as length spaces with the induced metric from $M_i$. Then the glued
space $M_1\cup_{i}M_2$ is an Alexandrov space with curvature $\geq
k$. }

\vskip1mm

We know that the Gluing Theorem is just the Doubling Theorem by Perel'man if $M_1=M_2$ ([Pe]).
If the glued space $M=M_1\cup_{i}M_2$ in the Gluing Theorem is a suspension, cone or join, we find that
there are some restrictions on the structures of $M_i$. In this paper, we make clear these restrictions.

%%%%%%%%%%%%%%%%%%%%%%%%%%%%%%%%%%%% Section 0 %%%%%%%%%%%%%%%%%%%%%%%%%%%%%%%%%%%%%%%%

\setcounter{section}{-1}

\section{Notations and main theorem}

We first make some conventions on the notations in the paper.

$\bullet$ Let $\mathcal{A}^n(k)$ denote the set of $n$-dimensional Alexandrov spaces with curvature $\geq k$.
In the paper, the spaces in $\mathcal{A}(k)$ are always assumed to be complete without special remark.

$\bullet$ Let $S(X)$, $C(X)$ and $X*Y$ denote the (spherical) suspension and cone over some $X\in\mathcal{A}(1)$, and the join between $X$ and $Y\in \mathcal{A}(1)$ respectively. They are defined as follows ([BGP]) (we will use $|xy|$ to denote the distance between $x$ and $y$ in $X$).

$\bullet$  $S(X)$ is the quotient space $X\times[0,\pi]/\sim$, where
$(x_1,a_1)\sim(x_2,a_2)\Leftrightarrow a_1=a_2=0$ or $a_1=a_2=\pi$,
with the metric
$$\cos|p_1p_2|=\cos a_1\cos a_2+\sin a_1\sin a_2\cos|x_1x_2|,$$ where $p_i=[(x_i,a_i)]$
(the class of $(x_i,a_i)$ in $X\times[0,\pi]/\sim$).

$\bullet$ $C(X)$ is the quotient space $X\times[0,+\infty)/\sim $, where $(x_1,a_1)\sim(x_2,a_2)\Leftrightarrow a_1=a_2=0$, with the metric
$$|p_1p_2|^2=a_1^2+a_2^2-2a_1a_2\cos|x_1x_2|,$$ where $p_i=[(x_i,a_i)]$.

$\bullet$ $X*Y$ is the quotient space $X\times Y\times[0,\frac\pi2]/\sim$, where $(x_1,y_1, a_1)\sim(x_2,y_2, a_2)\Leftrightarrow a_1=a_2=0$ and $x_1=x_2$ or $a_1=a_2=\frac\pi2$ and $y_1=y_2$, with the metric
$$\cos|p_1p_2|=\cos a_1\cos a_2\cos|x_1x_2|+\sin a_1\sin a_2\cos|y_1y_2|,$$ where $p_i=[(x_i,y_i,a_i)]$. Obviously,  $S(X)=\{p_1, p_2\}*X$ with $|p_1p_2|=\pi$.

Due to $X,Y\in \mathcal{A}(1)$, we have that $S(X), X*Y\in
\mathcal{A}(1)$ and $C(X)\in \mathcal{A}(0)$ ([BGP]). In Appendix of
the paper, we will discuss how to get the metric of $X*Y$, and
supply some basic properties of $X*Y$.

$\bullet$ Let $M_i\in \mathcal{A}(1)$ or $\mathcal{A}(0)$ ($i=1,2$)
with nonempty boundary $\partial M_i$. As in the Gluing Theorem, let
$\partial M_1$ be isometric to $\partial M_2$ (denoted by $\partial
M_1\stackrel{\rm{iso}}{\cong}\partial M_2$), and let
$M_1\cup_{\partial M_i}M_2$ denote the glued space. In the paper,
let $|\cdot|_i$ always denote the metric of $M_i$.

\vskip2mm

Now we formulate our main theorems as follows.

\vskip2mm

{\noindent \bf Theorem A}  {\it Let $M_i\in \mathcal{A}^n(1)$ with
boundaries $\partial M_1\stackrel{\rm{iso}}{\cong}\partial M_2$. If
there exist $p_i\in M_i$ such that $|p_1x|_1+|p_2x|_2\geq\pi$ for
any $x\in \partial M_i$, then one of the following cases holds:

\noindent{\rm(i)} $p_1,p_2\in \partial M_i$ and there exist $X_i\in \mathcal{A}(1)$ with boundaries
$\partial X_1\stackrel{\rm{iso}}{\cong}\partial X_2$ such that
$M_i=\{p_1,p_2\}*X_i$;

\noindent{\rm(ii)} $p_i\in M_i^\circ$\footnote{In the paper, we let $X^\circ$ denote the interior part of $X$.\hfill{$\,$}}
and $\partial M_i$ is convex\footnote{We say that $N$ is convex in $M$ if for any $x,
x'\in N$ there exists a geodesic $[xx']$ which falls in $N$. Of course, a convex subset
in $M\in \mathcal{A}(k)$ also belongs to $\mathcal{A}(k)$.} in $M_i$, and there exist $q_i\in M_i^\circ$ such that
$M_i=\{q_i\}*\partial M_i$, and thus $M_1\cup_{\partial
M_i}M_2=S(\partial M_i)$.}

\vskip1mm

Note that if $n=1$ in Theorem A, then $\partial M_i$ consists of two
points with distance equal to $\pi$; in this special case, we also
say that $\partial M_i$ is convex. It is easy to see that Theorem A
has the following two corollaries.

\vskip1mm

{\noindent \bf Corollary 0.1}  {\it Let $M_i\in \mathcal{A}^n(1)$
with boundaries $\partial M_1\stackrel{\rm{iso}}{\cong}\partial
M_2$. If $M_1\cup_{\partial M_i}M_2=S(X)$ for some
$X\in\mathcal{A}(1)$, then one of
the following cases holds:

\noindent{\rm(i)} there exist $X_i\in \mathcal{A}(1)$ with  $\partial
X_1\stackrel{\rm{iso}}{\cong}\partial X_2$ such that $M_i=S(X_i)$
and $X=X_1\cup_{\partial X_i}X_2$;

\noindent{\rm(ii)}  $\partial M_i$ is is convex in $M_i$ and  there
exist $q_i\in M_i^\circ$ such that $M_i=\{q_i\}*\partial M_i$, and
thus $X=\partial M_i$.}

\vskip1mm

{\noindent \bf Corollary 0.2}  {\it Let $M\in \mathcal{A}^n(1)$ with
nonempty boundary $\partial M$. If there exists $p\in M$ such that
$|px|\geq\frac{\pi}{2}$ for any $x\in\partial M$, then $\partial M$
is convex in $M$ and $M=\{p\}*\partial M$.}

\vskip1mm

Of course, Corollary 0.2 can also be derived by the Doubling
Theorem.

\vskip1mm

{\noindent \bf Remark 0.3} In the proof of Theorem A, the
quasi-geodesic ([PP],[Pet]) will be used. However, only in order to
prove Corollary 0.1 (and 0.2), we must not use the quasi-geodesic;
we can use the idea of proving Theorem B (which does not use the
quasi-geodesic) to give the proof.

\vskip2mm

{\noindent \bf Theorem B} {\it Let $M_i\in \mathcal{A}^n(0)$ with
boundaries $\partial M_1\stackrel{\rm{iso}}{\cong}\partial M_2$. If
$M_1\cup_{\partial M_i}M_2=C(X)$ with vertex $O$, where
$X\in\mathcal{A}^{n-1}(1)$ without boundary, then one of the following
holds:

\noindent{\rm(i)} there exists $X'\in\mathcal{A}(1)$ on which there is an isometrical
$\Bbb Z_2$-action that naturally induces an isometrical
$\Bbb Z_2$-action on $S(X')$ such that $X=S(X')/\Bbb Z_2$;

\noindent{\rm(ii)} there exist $X_1, X'\in \mathcal{A}(1)$ with $\partial
X_1\stackrel{\rm{iso}}{\cong}X'$ such that
$X=X_1\cup_{X'}(\{\xi\}*X')$;

\noindent{\rm(iii)} there exist $X_i\in \mathcal{A}(1)$ with $\partial
X_1\stackrel{\rm{iso}}{\cong}\partial X_2$ such that
$X=X_1\cup_{\partial X_i}X_2$.

\noindent Moreover, in the former two cases $O\in M_1^\circ$ or $M_2^\circ$, say $M_1^\circ$,
and then $M_2$ can be isometrically embedded into $C(\{\xi\}*X')$
and $\partial M_2$ is parallel to $\partial (C(\{\xi\}*X'))$; in the third case, $O\in \partial M_i$
and $M_i=C(X_i)$.}

\vskip1mm

Refer to Definition 2.5.2 and Corollary 2.5.3 for ``$\partial M_2$ is parallel to $\partial (C(\{\xi\}*X'))$''.

{\noindent\bf Example 0.4}  In Theorem B, suppose that $X$ is a
circle with perimeter $\ell$, then:

\noindent(0.4.1) If $\ell<\pi$, then only case (iii) occurs. And in
this case, $M_i$ are two sectors in $\Bbb R^2$ whose angles at the
vertex are $\alpha_i$ with $\alpha_1+\alpha_2=\ell$.

\noindent(0.4.2) If $\ell=\pi$, then case (i) or (iii) occurs. In
case (iii), $M_i$ are the same as in (0.4.1). In case (i), $X'$
consists of two points with distance $\le\pi$, so $\{\xi\}*X'$ is
just an arc of length $\pi$. Then $C(\{\xi\}*X')$ is just a half
plane; and thus $M_2$ is also a half plane in $C(\{\xi\}*X') $ with
$\partial M_2$ being a straight line which is parallel to $\partial
(C(\{\xi\}*X'))$.

\noindent(0.4.3) If $\ell>\pi$ (of course $\ell\leq2\pi$), then case
(ii) or (iii) occurs. In case (iii), $M_i$ are the same as in
(0.4.1), but it should be added that $\alpha_i\leq\pi$. In case
(ii), $X'$ and $M_2$ are the same as in case (i) in (0.4.2), and
$X_1$ is an arc of length $\ell-\pi$.

\vskip1mm

{\noindent \bf Remark 0.5} Given any $X\in \mathcal{A}(1)$ without
boundary satisfying (i)-(iii) in Theorem B, by the proof of Theorem
B (see Remark 2.5.1) we conclude that $C(X)$ can be split to two
parts $M_1, M_2\in \mathcal{A}(0)$ as in Theorem B. For example, let
$X$ be a circle with perimeter $\pi\leq\ell\leq 2\pi$. Select any
$\xi\in X$ (then $X'=\{\xi_1,\xi_2\}$ with $|\xi\xi_i|=\frac\pi2$),
and any point $p\in\gamma\setminus\{O\}$, where $\gamma\subset C(X)$
is the ray starting from the vertex $O$ with direction $\xi$. At
$p$, there is a unique local geodesic $\beta$ perpendicular to
$\gamma$. Note that $C(X)$ is split to $M_1$ and $M_2$ along
$\beta$, and $M_i$ with the induced intrinsic metric from $C(X)$
belongs to $\mathcal{A}^2(0)$. Moreover, $O\in M_1^\circ$ or
$M_2^\circ$, say $M_1^\circ$, and so $M_2$ is a half plane and
$\partial M_2=\beta$ is a straight line in $M_2$. However, if
$\ell<\pi$, then $\beta$ crosses itself, i.e. such splitting method
fails; in order to ensure that $M_i$ belongs to $\mathcal{A}^2(0)$,
according to Theorem B, $C(X)$ has to be split along two rays
starting from $O$ (i.e. $M_i$ are two sectors).

\vskip1mm

{\noindent \bf Remark 0.6} From Theorem A and the proof of Theorem B
(see Lemma 2.3), we conclude that $\partial M_i\in \mathcal{A}(1)$
in case (ii) of Theorem A (and Corollary 0.1), and $\partial M_i\in
\mathcal{A}(0)$ in cases (i) and (ii) of Theorem B. In general, we
do not know whether $\partial X$ is an Alexandrov space with
curvature bounded below for any $X\in \mathcal{A}(k)$. This is a BIG
conjecture (the Boundary Conjecture) in Alexandrov geometry ([BGP]):

\noindent{\it The boundary of a complete Alexandrov space with curvature $\geq k$
is a complete Alexandrov space with curvature $\geq k$ with respect to the induced intrinsic metric.}

\vskip2mm

{\noindent \bf Theorem C} {\it Let $M_i\in \mathcal{A}^n(1)$  with boundaries $\partial M_1\stackrel{\rm{iso}}{\cong}\partial M_2$.
If $M_1\cup_{\partial M_i}M_2=Y_1*Y_2$ for some $Y_1,Y_2\in\mathcal{A}(1)$ with empty boundary and $\diam(Y_2)<\pi$,
then there exist $X_i\in \mathcal{A}(1)$ with boundaries $\partial X_1\stackrel{\rm{iso}}{\cong}\partial X_2$ such that either
$Y_1=X_1\cup_{\partial X_i}X_2$ and $M_i=X_i*Y_2$, or $Y_2=X_1\cup_{\partial X_i}X_2$ and $M_i=Y_1*X_i$.}

\vskip1mm

{\noindent \bf Remark 0.7} Why do we add the condition
``$\diam(Y_2)<\pi$'' in Theorem C? Note that if
$\diam(Y_1)=\diam(Y_2)=\pi$, then there are
$Z_1,Z_2\in\mathcal{A}(1)$ with empty boundary such that
$Y_i=\{p_i,q_i\}*Z_i$ with $|p_iq_i|=\pi$, and thus $Y_1*Y_2=\Bbb
S^1*(Z_1*Z_2)$ where $\Bbb S^1$ has diameter equal to $\pi$. It then
follows that if $Y_1*Y_2$ can not be written to $\bar Y_1*\bar Y_2$
with $\diam(\bar Y_1)<\pi$ or $\diam(\bar Y_2)<\pi$, then
$Y_1*Y_2=\Bbb S^n$ with diameter equal to $\pi$, which implies that
$Y_j=\Bbb S^{n_j}$ with $n_1+n_2=n-1$. In such case, $M_1$ and $M_2$
have to be half spheres (by Corollary 0.1); however we can only say
that $M_i$ is isometric to $X_i*Y_2$ or $Y_1*X_i$ instead of that
$M_i=X_i*Y_2$ or $Y_1*X_i$. For example, let $M_1\cup_{\partial
M_i}M_2=\{p_1,p_2\}*\{q_1,q_2\}=\Bbb S^1$ with diameter equal to
$\pi$. Then $M_1$ can be any half of $\Bbb S^1$ (maybe not
$\{p_1,p_2\}*\{q_i\}$ or $\{p_i\}*\{q_1,q_2\}$).

\vskip2mm

Next we will formulate our mail tool---the Toponogov Comparison
Theorem (the essential geometry in Alexandrov spaces with curvature
bounded below).

$\bullet$ Let $[xy]$ denote a geodesic (i.e. shortest path) between $x$
and $y$ in $X\in\mathcal{A}(k)$.

$\bullet$ Let $\triangle pqr$ denote a triangle  in $X\in\mathcal{A}(k)$ consisting of three geodesics $[pq], [qr]$ and $[rp]$;
and let $\triangle\tilde p\tilde q\tilde r$ in $\Bbb S^2_k$ (the
complete and simply-connected 2-manifold of constant curvature $k$)
be a comparison triangle of $\triangle pqr$, i.e.  $|\tilde p\tilde
q|=|pq|,|\tilde p\tilde r|=|pr|$ and $|\tilde r\tilde q|=|rq|$ (recall that $|pq|+|pr|+|qr|\leqslant
2\pi/\sqrt k$ if $k>0$ ([BGP])).

$\bullet$ Let $p\prec^q_r$ denote a hinge in $X\in\mathcal{A}(k)$ consisting of two geodesics $[qp]$ and $[pr]$;
and let $\tilde p\prec^{\tilde q}_{\tilde r}$ in $\Bbb S^2_k$ be its comparison hinge, i.e. $|\tilde p\tilde
q|=|pq|$, $|\tilde p\tilde r|=|pr|$ and $\angle\tilde q\tilde p\tilde r=\angle qpr$.

\vskip2mm

\noindent {\bf Toponogov Comparison Theorem (TCT)}
{\it For any $\triangle pqr\subset X\in\mathcal{A}(k)$, we have $|ps|\geq|\tilde p\tilde s|,$ where $s\in[qr]$ and  $\tilde s\in[\tilde q\tilde r]\subset\triangle\tilde p\tilde q\tilde r$ with $|qs|=|\tilde q\tilde s|$. }

\vskip1mm

TCT has the following two equivalent versions:

\vskip1mm

\noindent {\bf TCT$'$} {\it For any $\triangle pqr\subset X$, we have  $\angle pqr\geq\angle \tilde p\tilde q\tilde r,$ $\angle qrp\geq\angle \tilde q\tilde r\tilde p$ and $\angle rpq\geq\angle \tilde r\tilde p\tilde q$.}

\vskip1mm

\noindent {\bf TCT$''$} {\it For any hinge $p\prec^q_r\subset X$ and its comparison hinge $\tilde p\prec^{\tilde q}_{\tilde r}$, we have  $|\tilde q\tilde r|\geq|qr|$.}

\vskip1mm

\noindent {\bf TCT for ``$=$'' ([GM])} {\it {\rm(i)} In TCT, if
there is $s\in[qr]^\circ$ such that $|ps|=|\tilde p\tilde s|,$ then
for any given geodesic $[ps]$ there exist unique two geodesics
$[pq]'$ and $[pr]'$ (maybe not $[pq]$ and $[pr]$) such that the
triangle formed by $[pq]'$, $[pr]'$ and $[qr]$ is isometric to its
comparison triangle.

\noindent{\rm(ii)} In TCT$'$ (resp. TCT$''$), if $\angle rpq=\angle
\tilde r\tilde p\tilde q$ (resp. $|\tilde q\tilde r|=|qr|$), then
there exists geodesic $[qr]'$ (maybe not $[qr]$  such that the
triangle formed by $[pq]$, $[pr]$ and $[qr]'$ is isometric to its
comparison triangle.}

\vskip1mm

On $X\in\mathcal{A}(k)$, a very important class of curves is the quasigeodesic ([PP], [Pet]). In particular, a local geodesic is a quasigeodesic.
We will use $x^\smallsmile y$ to denote a quasigeodesic between $x$ and $y$ in $X$. It is interesting that the following TCT still holds.

\noindent {\bf TCT with quasigeodesics ([PP])} {\it Let a geodesic $[pq]$ and quasigeodesic $p^\smallsmile r$ (resp. $[\tilde p\tilde q]$ and
a local geodesic $\tilde p^\smallsmile \tilde r$) form an angle equal to $\alpha$ at $p$ on $X\in\mathcal{A}(k)$ (resp. at $\tilde p$ on $\Bbb S^2_k$) with
$|\tilde p\tilde q|=|pq|$ and the length $\ell(\tilde p^\smallsmile \tilde r)=\ell(p^\smallsmile r)$. Then we have $|\tilde q\tilde r|\geq|qr|$.}

\vskip2mm

We will end this section with some other conventions (ref. [BGP]).

$\bullet$ Let $\Sigma_xX$ denote the direction space at $x\in X\in\mathcal{A}^{n}(k)$ which belongs to $\mathcal{A}^{n-1}(1)$.

$\bullet$ Let $\uparrow_x^y\in\Sigma_xX$ denote the direction of a geodesic $[xy]$ at $x\in X\in\mathcal{A}(k)$.

$\bullet$ We also use $|\cdot|_i$ to denote the metric of $\Sigma_xM_i$ if $x\in\partial M_i$.

$\bullet$ For convenience, we always use $N$ to denote the boundary
$\partial M_i$.

%%%%%%%%%%%%%%%%%%%%%%%%%%%%%%%%%%%% Section 1 %%%%%%%%%%%%%%%%%%%%%%%%%%%%%%%%%%%%%%%%

\section{Proof of Theorem A}

It is clear that Theorem A follows from Lemmas 1.1, 1,3 and 1.4
below. We will apply the induction on $\dim(M_i)$ (in Lemmas 1.2 and
1.3) to give the proof. Obviously, Theorem A is true if
$\dim(M_i)=1$, so we assume $\dim(M_i)>1$ in the rest of this
section.

\vskip2mm

{\noindent \bf 1.1\ \ On the case that one of $p_i$ belongs to $N$ ($=\partial M_i$)}

\vskip2mm

{\noindent \bf Lemma 1.1} {\it In Theorem A, if one of $p_i$ belongs to $N$, then both $p_1$ and $p_2$ belong to $N$; moreover, there exist $X_i\in \mathcal{A}^{n-1}(1)$ with boundaries $\partial X_1\stackrel{\rm{iso}}{\cong}\partial X_2$ such that $M_i=S(X_i)=\{p_1,p_2\}*X_i$.}

\vskip1mm

\noindent{\it Proof}. Without loss of generality, we assume that
$p_1\in N$.  Since $|p_1x|_1+|p_2x|_2\geq\pi$ for any $x\in N$, we
have $|p_2p_1|_2=\pi$. Hence, there exists $X_2\in
\mathcal{A}^{n-1}(1)$ with nonempty boundary such that
$M_2=\{p_1,p_2\}*X_2=S(X_2)$ and $\partial M_2=S(\partial
X_2)=\{p_1,p_2\}*\partial X_2$ (Corollary A.4.1 in Appendix).
Obviously, this implies that $p_1,p_2\in N\subset M_1$. Similarly,
$|p_1p_2|_1=\pi$, and there exists $X_1\in \mathcal{A}^{n-1}(1)$
such that $M_1=\{p_1,p_2\}*X_1=S(X_1)$ and $\partial
M_1=\{p_1,p_2\}*\partial X_1=S(\partial X_1)$. Note that $\partial
X_1\stackrel{\rm{iso}}{\cong}\partial X_2$ because $\partial
M_1\stackrel{\rm{iso}}{\cong}\partial M_2$. \hfill $\Box$

\vskip2mm

Due to Lemma 1.1, in the rest we only need to discuss the case that $p_i\in M_i^\circ$ for $i=1$ and 2.
Note that $p_i\in M_i^\circ$ implies that any $\uparrow_x^{p_i}\in(\Sigma_xM_i)^\circ$ for $x\in N$ ([BGP]).

\vskip2mm

{\noindent \bf 1.2\ \ A key observation: $|p_1x|_1+|p_2x|_2=\pi$}

\vskip2mm

{\noindent \bf Lemma 1.2} {\it In Theorem A, if in addition $p_i\in M_i^\circ$ for $i=1$ and $2$, then $|p_1x|_1+|p_2x|_2=\pi$ for any $x\in N$, and thus $|\uparrow_x^{p_1}\xi|_1+|\uparrow_x^{p_2}\xi|_2=\pi$ for any $\xi\in\Sigma_xN$.}

\vskip1mm

\noindent{\it Proof}.  The proof is inspired by [Pet].

Since $M_i$ is compact and $N$ is closed in $M_i$ ([BGP]), $N$
consists of finite components which are all compact. Let $N_0$ be
any component of $N$, and let $x_0\in N_0$  with
$$|p_1x_0|_1+|p_2x_0|_2=\min_{x\in N_0}\{|p_1x|_1+|p_2x|_2\}.$$ By
the first variation formula ([BGP]), we  have
$$|\uparrow_{x_0}^{p_1}\xi|_1+|\uparrow_{x_0}^{p_2}\xi|_2\geq\pi
\text{ for any } \xi\in\Sigma_{x_0}N.$$ Recall that
$\Sigma_{x_0}(\partial M_i)=\partial(\Sigma_{x_0} M_i)$ ([Pe]); and
thus since $\partial M_1\stackrel{\rm{iso}}{\cong}\partial M_2$,
$\partial(\Sigma_{x_0}
M_1)\stackrel{\rm{iso}}{\cong}\partial(\Sigma_{x_0} M_2)$ (cf.
[Pet]). Then by the induction on $\Sigma_{x_0} M_1$ and
$\Sigma_{x_0} M_2$ (note that $\uparrow_{x_0}^{p_i}\in (\Sigma_{x_0}
M_i)^\circ$ ([BGP])), we have
$$|\uparrow_{x_0}^{p_1}\xi|_1+|\uparrow_{x_0}^{p_2}\xi|_2=\pi \text{
for any } \xi\in\Sigma_{x_0}N.\eqno{(1.1)}$$ Now for any $x\in N_0$,
we select a shortest path $[xx_0]_N$ on $N$ between $x$ and $x_0$
which is a quasigeodesic in $M_i$ ([PP]). Due to $(1.1)$,
$|\uparrow_{x_0}^{p_1}\xi_0|_1+|\uparrow_{x_0}^{p_2}\xi_0|_2=\pi$,
where $\xi_0$ is the direction of $[xx_0]_N$ ([Pet]). On the unit
sphere $\Bbb S^2_1$, we select geodesics $[\tilde x_0\tilde p_i]$
and a local geodesic $\tilde x_0^\smallsmile\tilde x$ such that
$|\tilde x_0\tilde p_1|=|x_0p_1|_1, |\tilde x_0\tilde
p_2|=|x_0p_2|_2$ and the length $\ell(\tilde x_0^\smallsmile\tilde
x)=\ell([x_0x]_N)$, and $\tilde x_0^\smallsmile\tilde x$ is
perpendicular to $[\tilde x_0\tilde p_i]$ at $\tilde x_0$ and
$\angle\tilde p_1\tilde x_0\tilde p_2=\pi$. According to TCT with
quasigeodesics, we have $|p_1x|_1\leq|\tilde p_1\tilde x|$ and
$|p_2x|_2\leq|\tilde p_2\tilde x|$, and thus
$$|\tilde p_1\tilde x_0|+|\tilde p_2\tilde
x_0|=|p_1x_0|_1+|p_2x_0|_2\leq|p_1x|_1+|p_2x|_2\leq|\tilde p_1\tilde
x|+|\tilde p_2\tilde x|.$$ However, since $\angle\tilde p_1\tilde
x_0\tilde p_2=\pi$ and $|\tilde p_1\tilde x_0|+|\tilde p_2\tilde
x_0|\geq\pi$ on $\Bbb S^2_1$, it is not hard to see that
$$|\tilde p_1\tilde x|+|\tilde p_2\tilde x|\leq |\tilde p_1\tilde x_0|+|\tilde p_2\tilde x_0|,$$ and the `$=$' holds if and only if $|\tilde p_1\tilde x_0|+|\tilde p_2\tilde x_0|=\pi$. Therefore, we conclude that $$|\tilde p_1\tilde x_0|+|\tilde p_2\tilde x_0|=|p_1x_0|_1+|p_2x_0|_2=|p_1x|_1+|p_2x|_2=|\tilde p_1\tilde x|+|\tilde p_2\tilde x|=\pi.$$
Due to the arbitrary of $N_0$ as a component of $N$, we have $$|p_1x|_1+|p_2x|_2=\pi \text{ for any } x\in N,$$ and as a result
$$|\uparrow_{x}^{p_1}\xi|_1+|\uparrow_{x}^{p_2}\xi|_2=\pi \text{ for any } \xi\in\Sigma_xN \  \text{(see (1.1))}.$$
\hfill $\Box$

\vskip2mm

{\noindent \bf 1.3\ \ To find $q_i$ in Theorem A}

\vskip2mm

{\noindent \bf Lemma 1.3} {\it In Theorem A, if in addition $p_i\in M_i^\circ$ for $i=1$ and $2$, then there exist $q_i\in M_i^\circ$ such that $|q_1x|_1=|q_2x|_2=\frac{\pi}{2}$ for any $x\in N$.}

\vskip1mm

\noindent{\it Proof}. According to Lemma 1.2, $|p_1x|_1+|p_2x|_2=\pi$ for any $x\in N$. Then there exist $x_1$ and $x_2$ in $N$ such that (note that $N$ is compact)
$$\begin{aligned}&|p_1x_1|_1=\min\{|p_1x|_1|x\in N\},\ |p_2x_1|_2=\max\{|p_2x|_2|x\in N\},\\
&|p_1x_2|_1=\max\{|p_1x|_1|x\in N\},\ |p_2x_2|_2=\min\{|p_2x|_2|x\in N\}.
\end{aligned}$$
By the first variation formula ([BGP]), for any
$\xi\in\Sigma_{x_1}N$ and $\eta\in\Sigma_{x_2}N$ we have
$$|\uparrow_{x_1}^{p_1}\xi|_1\geq\frac\pi2,\
|\uparrow_{x_2}^{p_2}\eta|_2\geq\frac\pi2.$$ Then by the induction
(of Theorem A, see Corollary 0.2),
$$\Sigma_{x_i}M_i=\{\uparrow_{x_i}^{p_i}\}*\Sigma_{x_i}N$$ for $i=1$
and 2. Applying Lemma 1.2 again, we have
$$|\uparrow_{x_1}^{p_1}\xi|_1=|\uparrow_{x_1}^{p_2}\xi|_2=|\uparrow_{x_2}^{p_1}\eta|_1=|\uparrow_{x_2}^{p_2}\eta|_2=\frac\pi2.\eqno{(1.2)}$$
Similarly, we get
$$\Sigma_{x_1}M_2=\{\uparrow_{x_1}^{p_2}\}*\Sigma_{x_1}N \text{ and
} \Sigma_{x_2}M_1=\{\uparrow_{x_2}^{p_1}\}*\Sigma_{x_2}N,$$ which
implies that $$\angle p_1x_1x_2\leq\frac\pi2,\ \angle
p_1x_2x_1\leq\frac\pi2$$ in any triangle $\triangle p_1x_1x_2$ (as a
triangle in $M_1$). Let $\triangle\tilde p_1\tilde x_1\tilde x_2$ be
the comparison triangle of $\triangle p_1x_1x_2$. By TCT$'$,
$$\angle\tilde p_1\tilde x_1\tilde x_2\leq\angle p_1x_1x_2\leq\frac\pi2,\ \angle\tilde p_1\tilde x_2\tilde x_1\leq\angle p_1x_2x_1\leq\frac\pi2.\eqno{(1.3)}$$
It then follows that $$|p_1x_1|_1+|p_1x_2|_1=|\tilde p_1\tilde
x_1|+|\tilde p_1\tilde x_2|\leq\pi.$$  Similarly,
$|p_2x_1|_2+|p_2x_2|_2\leq\pi.$ On the other hand, note that
$|p_1x_1|_1+|p_2x_1|_2+|p_1x_2|_1+|p_2x_2|_2=2\pi$ (due to Lemma
1.2), hence we have
$$|p_1x_1|_1+|p_1x_2|_1=|p_2x_1|_2+|p_2x_2|_2=\pi.\eqno{(1.4)}$$
Note that $(1.3)$ and $(1.4)$ (together with the proof of Lemma 1.1)
imply that one of the following two cases holds:

Case 1: $|x_1x_2|_1<\pi$, $|x_1x_2|_2<\pi$ and $|p_ix_j|_i=\frac{\pi}{2}$.

Case 2:  $|x_1x_2|_1=|x_1x_2|_2=\pi.$

\noindent In Case 1, we let $q_i=p_i$. In Case 2,
$[x_1p_i]_i\cup[p_ix_2]_i$ is a geodesic of length $\pi$ in $M_i$
([BGP]), denoted by $[x_1x_2]_i$, and we let $q_i$ be the middle
point of $[x_1x_2]_i$.

It remains to show that $|q_ix|_i=\frac{\pi}{2}$ for any $x\in
N\setminus\{x_1,x_2\}$. Obviously, due to the choice of $x_1$ and
$x_2$, $|q_ix|_i=\frac{\pi}{2}$ in Case 1. In Case 2,
$[x_1x]_i\cup[xx_2]_i$ is a geodesic of length $\pi$ in $M_i$
between $x_1$ and $x_2$, and thus $[x_jx]_i\subset N$ (see footnote
5) which implies that $\angle q_ix_jx=\frac\pi2$ (see $(1.2)$); and
triangles $\triangle q_ix_jx$ are isometric to their comparison
triangles $\triangle\tilde q_i\tilde x_j\tilde x$ (Remark A.1.1 in
Appendix). It therefore follows that $|q_ix|_i=\frac\pi2.$ \hfill
$\Box$

\vskip2mm

{\noindent \bf 1.4\ \ To get $M_i=\{q_i\}*N$}

\vskip2mm

{\noindent \bf Lemma 1.4} {\it Let $M\in\mathcal{A}^{n}(1)$ with
nonempty boundary $N$. If there is a compact subset $A\subset M$
such that $|Ax|=\frac{\pi}{2}$ for any $x\in N$, then $A$ consists
of one point $q$ and $N$ is convex in $M$. Moreover, $M=\{q\}*N$.}

\vskip2mm

It is not hard to see that Lemma 1.4 is a corollary of the Doubling Theorem by Perel'man.
Here we supply a proof by the induction.

\vskip2mm

\noindent{\it Proof}. Obviously, the lemma is true when $n=1$, so we assume that $n>1$.

It follows from the first variation formula ([BGP]) that
$$|A'\xi|=\frac\pi2\text{ for any } \xi\in \Sigma_xN, \eqno{(1.5)}$$
where $A'=\{\uparrow_x^{a}|a\in A\}$. Hence, by the induction on
$\Sigma_xM$, $A'$ consists of one point $\uparrow_x^{q}$, which
implies that $A=\{q\}$ and there is a unique geodesic between $q$
and $x$; moreover,
$$\Sigma_x M=\{\uparrow_x^{q}\}*\Sigma_xN.\eqno{(1.6)}$$ This
implies that $$|\uparrow_x^{q}\eta|<\frac\pi2\text{ for any }
\eta\in (\Sigma_xM)^\circ.\eqno{(1.7)}$$

Due to $(1.5)$ and $(1.7)$, $\angle qx_1x_2, \angle
qx_2x_1\leq\frac\pi2$ in any triangle $\triangle qx_1x_2$ with $x_1,
x_2\in N$. On the other hand, in the comparison triangle
$\triangle\tilde  q\tilde x_1\tilde x_2$ of $\triangle qx_1x_2$,
$\angle \tilde q\tilde  x_1\tilde x_2=\angle\tilde  q\tilde
x_2\tilde x_1=\frac\pi2$ because $|\tilde q\tilde x_1|=|\tilde
q\tilde x_2|=\frac\pi2$. According to TCT$'$, it has to hold that
$$\angle qx_1x_2=\angle qx_2x_1=\frac\pi2.$$ Due to $(1.5)$ and
$(1.7)$ again, $\uparrow_{x_1}^{x_2}\in \Sigma_xN$, and thus
$[x_1x_2]\subset N$\footnote{For any $x,y$ in $M\in \mathcal{A}(k)$
with nonempty boundary, either $[xy]^\circ\subset M^\circ$ or
$[xy]^\circ\subset\partial M$ ([BGP]); in particular, if $x\in
\partial M$, then $\uparrow_x^y$ is an inner (resp. boundary)
direction in $\Sigma_xM$ if and only if the former (resp. latter)
case occurs.}, i.e. $N$ is convex in $M$.

Since there is a unique geodesic between $q$ and $x$ in $N$ which is
convex in $M$, in order to prove `$M=\{q\}*N$' we only need to show
that for any $y\in M$ there exists $x_0\in N$ such that $y\in
[x_0q]$ (see Remark A.3.4 in Appendix). In fact, $x_0$ is just the
point such that $|yx_0|=\min_{x\in N}\{|yx|\}$. Note that the first
variation formula implies that
$|\uparrow_{x_0}^{y}\xi|\geq\frac\pi2\text{ for any geodesic }
[x_0y] \text{ and } \xi\in \Sigma_{x_0}N.$ It then follows from
$(1.6)$ that $\uparrow_{x_0}^{y}=\uparrow_{x_0}^{q}$, i.e.
$[x_0y]\subseteq[x_0q]$ or vice versa. However, if $[x_0q]\subsetneq
[x_0y]$, then $[yq]\cup[qx]$ is a geodesic for any $x\in N$ (because
$|qx|=|qx_0|$), a contradiction. Hence, we have
$[x_0y]\subseteq[x_0q]$. \hfill $\Box$

\vskip1mm

%%%%%%%%%%%%%%%%%%%%%%%%%%%%%%%%%%%% Section 2 %%%%%%%%%%%%%%%%%%%%%%%%%%%%%%%%%%%%%%%

\section{Proof of Theorem B}

We will prove Theorem B according to two cases: $O\in M_1^\circ$ and
$O\in N\ (=\partial M_i)$. Subsections 2.1-2.6 are on the former
case, and subsection 2.7 is on the latter case.

In this section, we let $M$ denote $C(X)=M_1\cup_{\partial M_i}M_2$,
and let $\gamma_v$ denote the ray in $M$ starting from $O$ with
direction $v\in X$; and all lemmas are under the conditions in
Theorem B.

\vskip2mm

{\noindent \bf 2.1\ \ A basic observation on $N$}

\vskip2mm

For any local geodesic $c:[0,1]\to X$, we let $S_c=\{\gamma_{c(t)}|t\in[0,1]\}$. With respect to the induced metric, $S_c$ is just a (Euclidean) sector with $O$ being the vertex,
and the vertex angle of $S_c$ is equal to the length of $c$.

\vskip1mm

{\noindent \bf Lemma 2.1}  {\it Suppose that $O\in M_1^\circ$, and
that $c:[0,1]\to X$ is a local geodesic in $X$. If
$\gamma_{c(t)}\cap N\neq\emptyset$ for any $t\in [0,1]$, then
$S_c\cap N$ is a geodesic in the sector $S_c$.}

\vskip1mm

{\noindent \bf Sublemma 2.1.1}  {\it Suppose that $O\in M_1^\circ$.
Then for any $x\in M_1$, the geodesic $[Ox]_{M_1}$ in $M_1$ is just
the geodesic $[Ox]_{M}$ in $M$. As a result, for any $v\in X$,
$\gamma_v$ contains at most one point in $N$.}

\vskip1mm

\noindent{\it Proof}. Since $O\in M_1^\circ$,
$[Ox]_{M_1}\setminus\{x\}$ belongs to $M_1^\circ$ (see footnote 5).
Then $[Ox]_{M_1}\setminus\{x\}$ is a local geodesic in $M$ starting
from $O$. This implies that $[Ox]_{M_1}\setminus\{x\}$ has to lie in
the ray $\gamma_{\uparrow_O^{x}}$, and thus $[Ox]_{M_1}$ lies in the
ray $\gamma_{\uparrow_O^{x}}$. Hence, $[Ox]_{M_1}$ is just the
geodesic $[Ox]_{M}$ in $M$. Note that this (together with
$[Ox]_M\setminus\{x\}\subset M_1^\circ$) implies that $\gamma_v\cap
N$ contains at most one point for any $v\in X$.\hfill $\Box$

\vskip1mm

{\noindent \bf Sublemma 2.1.2} {\it Suppose that $O\in M_1^\circ$.
Then $|xy|_1=|xy|_M$ for any $x,y\in M_1$.}

\vskip1mm

\noindent{\it Proof}. Note that it suffices to show that
$|xy|_1\leq|xy|_M$ because $|xy|_1\geq|xy|_M$. From the definition
of the metric of the cone, we know that any $\triangle Oxy$ as a
triangle in $M$ is isometric to its comparison triangle
$\triangle\tilde O\tilde x\tilde y$ (in $\Bbb R^2$). It is clear
that $\angle xOy=\angle\tilde x\tilde O\tilde y$. On the other hand,
let $\bar\triangle Oxy$ be the triangle in $M_1$ with vertices $O,
x$ and $y$, and let $\triangle\bar O\bar x\bar y$ be its comparison
triangle. By TCT$'$, $\bar\angle xOy\geq\angle\bar x\bar O\bar y$,
where $\bar\angle xOy$ is the angle at vertex $O$ in $\bar\triangle
Oxy$. According to Sublemma 2.1.1, geodesics $[Ox]_{M_1}$ and
$[Oy]_{M_1}$ in $M_1$ are just the geodesics $[Ox]_{M}$ and
$[Oy]_{M}$ in $M$ respectively. Since  $O\in M_1^\circ$, $\angle
xOy=\bar\angle xOy$, so $\angle\bar x\bar O\bar y\leq\angle\tilde
x\tilde O\tilde y$. It then follows that
$$\hskip5cm|xy|_1=|\bar x\bar y|\leq |\tilde x\tilde y|=|xy|_M.\hskip4.5cm\Box$$

\vskip1mm

{\noindent \bf Sublemma 2.1.3} {\it For any $x_1, x_2\in N$, if
there exists a point $p\in M_2^\circ$ such that $|px_j|_2=|px_j|_M$
and $\uparrow_O^p$ lies in a geodesic
$[\uparrow_O^{x_1}\uparrow_O^{x_2}]\subset X$, then
$|x_1x_2|_2=|x_1x_2|_M$.}

\vskip1mm

\noindent{\it Proof}. Similarly, we only need to show that
$|x_1x_2|_2\leq|x_1x_2|_M$. Since $|px_j|_2=|px_j|_M$, each geodesic
$[px_j]_{M_2}$ in $M_2$ is a geodesic $[px_j]_{M}$ in $M$. If
$\uparrow_O^p\in[\uparrow_O^{x_1}\uparrow_O^{x_2}]^\circ$, then
there is a unique geodesic between $\uparrow_O^p$ and
$\uparrow_O^{x_j}$, and thus there is a unique geodesic between $p$
and $x_j$ in $M$ (due to the definition of the metric of $C(X)$); if
$\uparrow_O^p=\uparrow_O^{x_1}$ or $\uparrow_O^{x_2}$, say
$\uparrow_O^{x_1}$, then $[px_1]_{M_2}$$(=[px_1]_M)$ lies in the ray
$\gamma_{\uparrow_O^{x_1}}$. In any case, geodesics
$[px_j]_{M_2}$$(=[px_j]_M)$ belong to a sector
$S_{[\gamma_{\uparrow_O^{x_1}}\gamma_{\uparrow_O^{x_2}}]}$. Let
$\triangle px_1x_2$ be the triangle (in $M$) which lies in
$S_{[\gamma_{\uparrow_O^{x_1}}\gamma_{\uparrow_O^{x_2}}]}$. Note
that $\triangle px_1x_2$ itself is its comparison triangle (in $\Bbb
R^2$). On the other hand, let $\bar\triangle px_1x_2$ be a triangle
in $M_2$ containing the sides $[px_j]_{M_2}$ ($=[px_j]_M$). Since
$p\in M_2^\circ$, $\angle x_1px_2=\bar\angle x_1px_2$, where
$\bar\angle x_1px_2$ is the angle at vertex $p$ in $\bar\triangle
px_1x_2$. Let $\triangle\bar p\bar x_1\bar x_2$ be the comparison
triangle of $\bar\triangle px_1x_2$. By TCT$'$, $\bar\angle
x_1px_2\geq\angle\bar x_1\bar p\bar x_2$. Hence, $\angle\bar x_1\bar
p\bar x_2\leq \angle x_1px_2$, so we have
$$\hskip5cm|x_1x_2|_2=|\bar x_1\bar x_2|\leq |x_1x_2|_M.\hskip5cm\Box$$

\vskip1mm

{\noindent \bf Sublemma 2.1.4} {\it Assume that $O\in M_1^\circ$,
$p\in M_2^\circ$ and $[Op]_M\cap N=\{x_0\}$. Then  for any $x\in
N\setminus\{x_0\}$, there exists $x_1\in[px]_M\cap
N\setminus\{x_0\}$ such that $|px_1|_2=|px_1|_M$.}

\vskip1mm

\noindent{\it Proof}. Due to Sublemma 2.1.1,
$\gamma_{\uparrow_O^p}\cap N=\{x_0\}$, so
$x\not\in\gamma_{\uparrow_O^p}$. This implies that
$\gamma_{\uparrow_O^p}\cap[px]_M=\{p\}$. Then the sublemma follows.
\hfill $\Box$

\vskip2mm

\noindent{\it Proof of Lemma 2.1}.

Since $c$ is a local geodesic in $X$, for any $0<t_0<1$, there are
$0<t_1<t_2<1$ such that $t_0\in(t_1, t_2)$ and $c|_{[t_1,t_2]}$ is a
geodesic. In addition, we may assume that $c|_{[t_1,t_2]}$ is the
unique geodesic between $c(t_1)$ and $c(t_2)$. Note that
$\gamma_{c(t_0)}\subset S_{c|_{[t_1,t_2]}}\subset S_c$.

By Sublemma 2.1.1, $\gamma_{c(t)}$ contains only one point in $N$. Let $x_j=\gamma_{c(t_j)}\cap N$ for $j=0, 1$ and 2, and let $z$
be a point in $\gamma_{c(t_0)}\setminus[Ox_0]$ which belongs to $M_2^\circ$. Note that geodesic $[zx_j]_M$ lies in the sector
$S_{c|_{[t_0,t_j]}}\subset S_{c|_{[t_1,t_2]}}$. Then according to Sublemma 2.1.4, $t_j$ can be selected originally such that $|zx_j|_2=|zx_j|_M$ for $j=1$ and 2.
It therefore follows from Sublemmas 2.1.2 and 2.1.3 that
$$|x_1x_2|_1=|x_1x_2|_2=|x_1x_2|_M,$$
which implies that any geodesic $[x_1x_2]_{M_i}$ is a geodesic in $M$. Note that
there is a unique geodesic $[x_1x_2]_M$ between $x_1$ and $x_2$ in $M$ because $c|_{[t_1,t_2]}$ is the
unique geodesic between $c(t_1)$ and $c(t_2)$.
This implies that $$[x_1x_2]_{M_1}=[x_1x_2]_{M_2}=[x_1x_2]_M\subset N.$$
Note that $[x_1x_2]_M$ is just the geodesic between $x_1$ and $x_2$ in $S_{c|_{[t_1,t_2]}}$, so $$[x_1x_2]_M=N\cap S_{c|_{[t_1,t_2]}}.$$

Due to the arbitrary of $t_0\in(0,1)$ and the closeness of $N$, we
conclude that $N\cap S_{c}$ is a geodesic in $S_{c}$ with end points
belonging to $\gamma_{c(0)}$ and $\gamma_{c(1)}$. \hfill $\Box$

\vskip2mm

{\noindent \bf 2.2\ \ The non-compactness of $N$}

\vskip2mm

{\noindent \bf Proposition 2.2} {\it Suppose that $O\in M_1^\circ$. Then $N$ is not compact.}

\vskip1mm

{\noindent \bf Sublemma 2.2.1} {\it Suppose that $O\in M_1^\circ$.
Then $X_0:=\{\uparrow_O^x|x\in N\}$ is open in $X$.}

\vskip1mm

\noindent{\it Proof}.  If $X_0$ is not open in $X$, then there
exists $\uparrow_O^{x_0}\in X_0$ with a sequence $v_k\in X$
converging to $\uparrow_O^{x_0}$ such that the rays $\gamma_{v_k}$
(in $M$) contain no point in $N$. Note that $\gamma_{v_k}\subset
M_1$ because $O\in M_1$, and  $\gamma_{v_k}$ converges to the ray
$\gamma_{\uparrow_O^{x_0}}$ (in $M$). It follows that
$\gamma_{\uparrow_O^{x_0}}$ belongs to $M_1$. On the other hand,
$\gamma_{\uparrow_O^{x_0}}$ contains a unique boundary point $x_0$
(see Sublemma 2.1.1). This is impossible (see footnote 5). \hfill
$\Box$

\vskip2mm

\noindent{\it Proof of Proposition 2.2}.

We will derive a contradiction by assuming that $N$ is compact.
``$N$ is compact'' implies that there is $x_0\in N$ such that
$|Ox_0|_1=\max\{|Ox|_1|x\in N\}$. Note that there is a unique
geodesic between $O$ and $x_0$ (by Sublemma 2.1.1). By the first
variation formula,
$$|\uparrow_{x_0}^O\xi|_1\leq\frac{\pi}{2} \text { for any
}\xi\in\Sigma_{x_0}N.$$ {\bf Claim}: {\it In fact,
$$|\uparrow_{x_0}^O\xi|_1=\frac{\pi}{2}.\eqno{(2.1)}$$} Note that for any given
$z$ in $\gamma_{\uparrow_O^{x_0}}\setminus[Ox_0]$ ($\subset
M_2^\circ$ by Sublemma 2.1.1),
$$|Oz|=|Ox_0|_1+|x_0z|_2=\min\{|Ox|_1+|xz|_2|x\in N\},$$ so from the
first variation formula ([BGP])
$$|\uparrow_{x_0}^O\xi|_1+|\uparrow_{x_0}^z\xi|_2\geq\pi.\eqno{(2.2)}$$
On the other hand, $\uparrow_{x_0}^O\in(\Sigma_{x_0}M_1)^\circ$ and
$\uparrow_{x_0}^z\in(\Sigma_{x_0}M_2)^\circ$ (because $O\in
M_1^\circ$ and $z\in M_2^\circ$). It therefore follows from Theorem
A that both $\Sigma_{x_0}M_1$ and $\Sigma_{x_0}M_2$ are half
$S(\Sigma_{x_0}N)$. Note that (2.2) implies that
$|\uparrow_{x_0}^z\xi|_2\geq\frac{\pi}{2}$ because
$|\uparrow_{x_0}^O\xi|_1\leq\frac{\pi}{2}$. Then
$\Sigma_{x_0}M_2=\{\uparrow_{x_0}^z\}$ $*\Sigma_{x_0}N$, and thus
$|\uparrow_{x_0}^z\xi|_2=\frac{\pi}{2}$. Again by (2.2) and
`$|\uparrow_{x_0}^O\xi|_1\leq\frac{\pi}{2}$', we have
$|\uparrow_{x_0}^O\xi|_1=\frac{\pi}{2},$ i.e. the claim is verified.

Note that $X_0=\{\uparrow_O^x|x\in N\}$ is a closed subset in $X$
under the assumption ``$N$ is compact'', which together with
Sublemma 2.2.1 implies that $X_0=X$. Then any geodesic
$[\uparrow_O^{x_0}\uparrow_O^{x}]$ in $X$ for any $x\in
N\setminus\{x_0\}$ satisfies Lemma 2.1, so
$S_{[\uparrow_O^{x_0}\uparrow_O^{x}]}\cap N$ is the geodesic
$[x_0x]$ between $x_0$ and $x$ in
$S_{[\uparrow_O^{x_0}\uparrow_O^{x}]}$. Due to (2.1), $[x_0x]$ is
perpendicular to $\gamma_{\uparrow_O^{x_0}}$ at $x_0$ in
$S_{[\uparrow_O^{x_0}\uparrow_O^{x}]}$. It then follows that
$$|Ox|=\sqrt{|Ox_0|^2+|xx_0|^2}>|Ox_0|,$$ which contradicts the choice of $x_0$. \hfill $\Box$

\vskip2mm

{\noindent \bf 2.3\ \ The convexity of $N$}

\vskip2mm

{\noindent \bf Lemma 2.3} {\it Suppose that $O\in M_1^\circ$. Then

\noindent{\rm (i)} $N$ is totally convex in $M_2$;

\noindent{\rm (ii)} $N$ is locally totally convex in $M_1$ and $M$.}

\noindent{\it Proof}. (i)\ \ We need show that any geodesic
$[xy]_{M_2}$ lies in $N$ for all $x, y\in N$. Note that if there is
an inner point of $[xy]_{M_2}$ belonging to $N$, then
$[xy]_{M_2}\subseteq N$ (see footnote 5). In the following, we will
derive a contradiction by assuming that $[xy]^\circ_{M_2}\subset
M_2^\circ$.

Observe that `$[xy]^\circ_{M_2}\subset M_2^\circ$' implies that any
geodesic $[x_ky_k]_{M_2}\subset[xy]^\circ_{M_2}$ is a local geodesic
in $M$, and thus $\sigma_k=\{\uparrow_O^z|z\in[x_ky_k]_{M_2}\}$ is a
local geodesic in $X$. Since $z\in M_2^\circ$ for any
$z\in[x_ky_k]_{M_2}$ and $O\in M_1^\circ$,
$\gamma_{\uparrow_O^z}\cap N\neq\emptyset$. It therefore follows
from Lemma 2.1 that $S_{\sigma_k}\cap N$ is a geodesic in the sector
$S_{\sigma_k}$, denoted by $[x_k'y_k']$. Now let $x_k$ and $y_k$
converge to $x$ and $y$ respectively. Note that $x'_k$ and $y'_k$
also converge to $x$ and $y$ respectively because $N$ is closed and
$\gamma_{\uparrow_O^x}\cap N=\{x\}$ and $\gamma_{\uparrow_O^y}\cap
N=\{y\}$ (see Sublemma 2.1.1). On the other hand, note that
$[x_ky_k]_{M_2}$ is also a geodesic in $S_{\sigma_k}$. Hence, the
two geodesics $[x_ky_k]_{M_2}$ and $[x_k'y_k']$ in $S_{\sigma_k}$
have to be the same one. This yields a contradiction because
$[x_ky_k]_{M_2}\subset M_2^\circ$ and $[x'_ky'_k]\subset N$.

\vskip1mm

(ii)\ \ From (i), $N$ with the induced metric from $M_2$ also
belongs to $\mathcal A(0)$. We need show that there is a
neighborhood $U_x\subset N$ of any $x\in N$ such that any geodesic
$[ab]_N$ is a geodesic in $M_1$ (so in $M$ due to (i)) for any
$a,b\in U_x$. Since $X_0$ is open in $X$ (see Sublemma 2.2.1), there
is a neighborhood $V$ of $\uparrow_O^{x}$ such that $[v_1v_2]\subset
X_0$ for any $v_1, v_2\in V$. Let $x_j=\gamma_{v_j}\cap N$ (see
Sublemma 2.1.1). By Lemma 2.1, $S_{[v_1v_2]}\cap N$ is a geodesic
between $x_1$ and $x_2$ in $S_{[v_1v_2]}$, so in $M$. Together with
Sublemma 2.1.2, this implies that
$$|x_1x_2|_N=|x_1x_2|_M=|x_1x_2|_1.$$ Hence, any $[x_1x_2]_N$ is a
geodesic in $M_1$, so it suffices to let $U_x=\{\gamma_v\cap N|v\in
V\}$. \hfill $\Box$

\vskip1mm

{\noindent \bf Remark 2.3.1}  In Lemma 2.3, $N$ may not be convex in
$M_1$ and $M$ (see the case `$\pi\leq\ell<2\pi$' of the example in
Remark 0.5). Note that, given a geodesic $[xy]_{M_1}$ in $M_1$ with
$x,y\in N$, for any $z\in[xy]_{M_1}^\circ$ the ray
$\gamma_{\uparrow_O^z}$ may contain no point in $N$; and thus the
argument in the proof of (i) of Lemma 2.3 fails when one try to
prove that `$N$ is convex in $M_1$'.

\vskip2mm

{\noindent \bf 2.4\ \ The nearest point to $O$ in $N$}

\vskip2mm

{\noindent \bf Lemma 2.4} {\it Suppose that $O\in M_1^\circ$. Then there exists a unique point $x_0\in N$ such that $|Ox_0|_1=\min\{|Ox|_1|x\in N\}$. Moreover,

\noindent{\rm (i)} $|\uparrow_O^{x_0}\uparrow_O^{x}|<\frac\pi2$ for any $x\in N$;

\noindent{\rm (ii)} any geodesic $[x_0x]_M$ lies in $N$ for any $x\in N$,
and thus $|x_0x|=|Ox_0|\tan|\uparrow_O^{x_0}\uparrow_O^{x}|$ and
$|Ox_0|=|Ox|\cos|\uparrow_O^{x_0}\uparrow_O^{x}|$;

\noindent{\rm (iii)} $N=\{p\in
M||\uparrow_{x_0}^{p}\uparrow_{x_0}^{O}|=\frac{\pi}{2} \text{ for
some } [x_0p]_M\}$.}

\noindent{\it Proof}. Since $O\in M_1^\circ$, there exists $x_0\in
N$ such that $|Ox_0|_1=\min\{|Ox|_1|x\in N\}$. By the first
variation formula, $|\uparrow_{x_0}^O\xi|_1\geq\frac{\pi}{2}$ for
any $\xi\in\Sigma_{x_0}N$, and thus
$\Sigma_{x_0}M_1=\{\uparrow_{x_0}^O\}$$*\Sigma_{x_0}N$ (see
Corollary 0.2). This implies that
$|\uparrow_{x_0}^O\xi|_1=\frac{\pi}{2}$ in fact. According to
Sublemma 2.1.2, the distance $|\eta_1\eta_2|_1$ for any
$\eta_1,\eta_2\in\Sigma_{x_0}M_1$ is just the distance
$|\eta_1\eta_2|$ in $\Sigma_{x_0}M$. It then follows that, in
$\Sigma_{x_0}M$, for any $\xi\in\Sigma_{x_0}N$
$$|\uparrow_{x_0}^O\xi|=\frac{\pi}{2}.\eqno{(2.3)}$$

(i)\ \ According to (ii) of Lemma 2.3, $[x_0x]_N$ is a local
geodesic in $M$, so $\sigma:=\{\uparrow_O^y|y\in [x_0x]_N\}$ is a
local geodesic in $X$. By Lemma 2.1, we conclude that $[x_0x]_N$ is
the geodesic between $x_0$ and $x$ in the sector $S_{\sigma}$. Since
$[x_0x]_N$ is perpendicular to $[x_0O]$ at $x_0$ (see (2.3)), the
length of $\sigma$ (i.e. the vertex angle of $S_\sigma$)
$\ell(\sigma)<\frac{\pi}{2}$. Of course,
$$|\uparrow_O^{x_0}\uparrow_O^{x}|\le\ell(\sigma)<\frac\pi2.$$

\vskip1mm

(ii)\ \ By Sublemma 2.2.1, there exists $y_1\in[x_0x]_M^\circ$ such
that $\{\uparrow_O^y|y\in[x_0y_1]_M\}\subset X_0$. Let
$x_1=N\cap\gamma_{\uparrow_O^{y_1}}$ (see Sublemma 2.1.1). Due to
Lemma 2.1, $S_{[\uparrow_O^{x_0}\uparrow_O^{x_1}]}\cap
N=[x_0x_1]_M$; and thus by (2.3)
$$|Ox_1|=\frac{|Ox_0|}{\cos|\uparrow_O^{x_0}\uparrow_O^{x_1}|}
\leq\frac{|Ox_0|}{\cos|\uparrow_O^{x_0}\uparrow_O^{x}|},$$ where
$|\uparrow_O^{x_0}\uparrow_O^{x}|<\frac\pi2$ (due to (i)). By
Sublemma 2.2.1 again, there exists $y_2\in[y_1x]_M^\circ$ such that
$\{\uparrow_O^y|y\in[x_0y_2]_M\}\subset X_0$. Similarly, let
$x_2=N\cap\gamma_{\uparrow_O^{y_2}}$, and we have
$S_{[\uparrow_O^{x_0}\uparrow_O^{x_2}]}\cap N=[x_0x_2]_M$ and
$$|Ox_2|=\frac{|Ox_0|}{\cos|\uparrow_O^{x_0}\uparrow_O^{x_2}|}\leq\frac{|Ox_0|}{\cos|\uparrow_O^{x_0}\uparrow_O^{x}|}.$$
Repeating the above process, we obtain that
$S_{[\uparrow_O^{x_0}\uparrow_O^{x}]}\cap N=[x_0x]_M$, and
$$|x_0x|=|Ox_0|\tan|\uparrow_O^{x_0}\uparrow_O^{x}|\text{ and }
|Ox_0|=|Ox|\cos|\uparrow_O^{x_0}\uparrow_O^{x}|.$$

(iii)  Since any $[x_0x]_M\subset N$ for any $x\in N$,
$N\subseteq\{p\in M|
|\uparrow_{x_0}^{p}\uparrow_{x_0}^{O}|=\frac{\pi}{2}\text{ for some
} [x_0p]_M\}$ due to (2.3). Then it suffices to show that if
$|\uparrow_{x_0}^{p}\uparrow_{x_0}^{O}|=\frac{\pi}{2}$ (in
$\Sigma_{x_0}M$) for some $[x_0p]_M$, then $p\in N$. We have proved
that $\Sigma_{x_0}M_1=\{\uparrow_{x_0}^O\}*\Sigma_{x_0}N$ in the
beginning of the proof. It also holds that
$\Sigma_{x_0}M_2=\{\uparrow_{x_0}^z\}*\Sigma_{x_0}N$, where $z\in
\gamma_{\uparrow_O^{x_0}}\setminus[Ox_0]$ (of course, $z\in
M_2^\circ$ due to Sublemma 2.1.1 and footnote 5). In fact, according
to the proof of (ii), it is not hard to see that
$|zx_0|=|zx|\cos|\uparrow_z^{x_0}\uparrow_z^{x}|$. This implies that
$|zx_0|_2=\min\{|zx|_2|x\in N\}$ (note that $|zx_0|_2=|zx_0|$ and
$|zx|_2\geq|zx|$), and thus similarly we have
$\Sigma_{x_0}M_2=\{\uparrow_{x_0}^z\}*\Sigma_{x_0}N$. Then
``$|\uparrow_{x_0}^{p}\uparrow_{x_0}^{O}|=\frac{\pi}{2}$'' implies
that $\uparrow_{x_0}^p\in\Sigma_{x_0}N$. Now we assume that
$p\not\in N$, i.e. $p\in M_i^\circ$ for $i=1$ or 2. Let $q$ be the
nearest point in $N$ to $p$ along the geodesic $[px_0]_M$. Due to
(ii), $[x_0q]_M\subset N$, so $[x_0p]_M\subset M_i$. This implies
that $\uparrow_{x_0}^p\in(\Sigma_{x_0}M_i)^\circ$ (see footnote 5),
which contradicts ``$\uparrow_{x_0}^p\in\Sigma_{x_0}N$''.

\vskip1mm

At last we show the uniqueness of $x_0$. Let $x'_0\in N$ be another point such that $|Ox'_0|_1=\min\{|Ox|_1|x\in N\}$.
Due to (ii) and (iii), any geodesic $[x_0x'_0]_M$ lies in $N$, and $|\uparrow_{x_0}^O\uparrow_{x_0}^{x'_0}|=|\uparrow_{x'_0}^{x_0}\uparrow_{x'_0}^O|=\frac{\pi}{2}.$
This is impossible because $[Ox_0]_M$, $[Ox'_0]_M$ and $[x_0x'_0]_M$ lie in the sector $S_{[\uparrow_O^{x_0}\uparrow_{O}^{x'_0}]}$.
\hfill $\Box$

\vskip2mm

{\noindent \bf 2.5\ \ The structures of
$X_0:=\{\uparrow_{O}^{x}|x\in N\}$ and $M_2$}

\vskip2mm

{\noindent \bf Lemma 2.5} {\it Suppose that $O\in M_1^\circ$, and
that $x_0\in N$ with $|Ox_0|_1=\min\{|Ox|_1|x\in N\}$. Then $X_0$
with the induced metric from $X$ is isometric to
$(\{\uparrow_{O}^{x_0}\}*\Sigma_{\uparrow_{O}^{x_0}}X)^\circ$.}

\vskip1mm

\noindent{\it Proof}. Since $N$ is locally totally convex in $M$
(Lemma 2.3), $X_0$ is locally totally convex in $X$. Hence, $X_0$
with the induced metric from $X$ belongs to $\mathcal{A}(1)$, and is
of dimension $n-1$ by Sublemma 2.2.1. Moreover, for any geodesic $[x_1x_2]_N$ in $N$,
$\{\uparrow_O^x|x\in[x_1x_2]_N\}$ is a geodesic between
$\uparrow_{O}^{x_1}$ and $\uparrow_{O}^{x_2}$ in $X_0$.
(Here, we note that $N$ is complete because $M_i$ is complete and $N$ is closed and convex in $M_2$.
However $X_0$ is not complete, otherwise $N$ will be compact which contradicts
Proposition 2.2).

Denote by $|\cdot|_0$ the induced metric on $X_0$
from $X$. We first note that
$|\uparrow_{O}^{x_0}\uparrow_{O}^{x}|_0=|\uparrow_{O}^{x_0}\uparrow_{O}^{x}|$
for any $x\in N$ because any geodesic $[x_0x]_N$ is a geodesic in $M$ ((ii) of Lemma 2.4),
then by (i) of Lemma 2.4 we have
$$|\uparrow_{O}^{x_0}\uparrow_{O}^{x}|_0<\frac{\pi}{2}.$$ For any
$x_1, x_2\in N\setminus\{x_0\}$,
$|\uparrow_{O}^{x_1}\uparrow_{O}^{x_2}|_0$ is the length of the
local geodesic $\{\uparrow_O^x|x\in[x_1x_2]_N\}$ in $X$ for any geodesic $[x_1x_2]_N$. {\bf Claim
1}: {\it $$\cos|\uparrow_{O}^{x_1}\uparrow_{O}^{x_2}|_0=
\cos|\uparrow_{O}^{x_0}\uparrow_{O}^{x_1}|_0\cos|\uparrow_{O}^{x_0}\uparrow_{O}^{x_2}|_0+
\sin|\uparrow_{O}^{x_0}\uparrow_{O}^{x_1}|_0\sin|\uparrow_{O}^{x_0}\uparrow_{O}^{x_2}|_0\cos\alpha,\eqno{(2.4)}$$
where $\alpha$ is the angle between
$[\uparrow_{O}^{x_0}\uparrow_{O}^{x_1}]$ and
$[\uparrow_{O}^{x_0}\uparrow_{O}^{x_2}]$ in $X_0\subset X$.}

In order to verify Claim 1, it suffices to prove that the equality
holds if $\uparrow_{O}^{x_j}$ is replaced by
$\uparrow_O^{x_j'}\in[\uparrow_{O}^{x_0}\uparrow_{O}^{x_j}]^\circ$
($j=1, 2$) in (2.4) (then we only need to let $\uparrow_O^{x_j'}$
converge to $\uparrow_O^{x_j}$). Note that there is a unique
geodesic between $\uparrow_{O}^{x_0}$ and $\uparrow_{O}^{x_j'}$. For
convenience, we still use $x_j$ to denote $x_j'$. Let
$\triangle\tilde\uparrow_{O}^{x_0}\tilde\uparrow_{O}^{x_1}\tilde\uparrow_{O}^{x_2}$
be the comparison triangle (on the unit sphere) of
$\triangle\uparrow_{O}^{x_0}\uparrow_{O}^{x_1}\uparrow_{O}^{x_2}$ in
$X_0$. Since we now have the condition that there is a unique
geodesic between $\uparrow_{O}^{x_0}$ and $\uparrow_{O}^{x_j}$,
according to TCT for `$=$', in order to prove (2.4) we only need to
show that
$$|\uparrow_{O}^{x_0}\uparrow_{O}^{y}|_0=|\tilde\uparrow_{O}^{x_0}\tilde\uparrow_{O}^{y}|,\eqno{(2.5)}$$
where $y\in [x_1x_2]_N^\circ$ and
$\tilde\uparrow_{O}^{y}\in[\tilde\uparrow_{O}^{x_1}\tilde\uparrow_{O}^{x_2}]$
with
$|\tilde\uparrow_{O}^{x_1}\tilde\uparrow_{O}^{y}|=|\uparrow_{O}^{x_1}\uparrow_{O}^{y}|_0$.
Due to Lemma 2.4,
$$|Ox_0|=|Oy|\cos\angle x_0Oy=|Oy|\cos|\uparrow_{O}^{x_0}\uparrow_{O}^{y}|_0.\eqno{(2.6)}$$
On the other hand, we consider the cone
$C(\triangle\tilde\uparrow_{O}^{x_0}\tilde\uparrow_{O}^{x_1}\tilde\uparrow_{O}^{x_2})$
with vertex $\tilde O$, a part of the 3-dimensional Euclidean space.
In this cone, we select $\tilde x_j$ with $j=0,1,2$ such that
$\uparrow_{\tilde O}^{\tilde x_j}=\tilde\uparrow_{O}^{x_j}$ and
$|\tilde O\tilde x_j|=|Ox_j|$. Note that for $\tilde y\in [\tilde
x_1\tilde x_2]$ with $\uparrow_{\tilde O}^{\tilde
y}=\tilde\uparrow_{O}^{y}$,
$$|Ox_0|=|\tilde O\tilde x_0|=|\tilde O\tilde y|\cos|\uparrow_{\tilde O}^{\tilde x_0}\uparrow_{\tilde O}^{\tilde y}|=
|\tilde O\tilde y|\cos|\tilde \uparrow_{O}^{x_0}\tilde
\uparrow_{O}^{y}|.\eqno{(2.7)}$$ Note that $\triangle Ox_1x_2$ (a
triangle in the sector $S_{\{\uparrow_O^x|x\in[x_1x_2]_N\}}$) is
isometric to the Euclidean triangle $\triangle\tilde  O\tilde
x_1\tilde x_2$, so $|Oy|=|\tilde O\tilde y|$. Then (2.5) follows
from (2.6) and (2.7), so Claim 1 holds.

Note that Claim 1 enables us to construct an isometrical embedding
$$i:X_0\longrightarrow\{\uparrow_{O}^{x_0}\}*\Sigma_{\uparrow_{O}^{x_0}}X \text{ defined by }
\uparrow_O^{x}\longmapsto[(\uparrow_{O}^{x_0},\uparrow_{\uparrow_{O}^{x_0}}^{\uparrow_{O}^{x}},|\uparrow_{O}^{x_0}\uparrow_{O}^{x}|_0)].$$

\noindent{\bf Claim 2}: {\it Both $i(X_0)$ and $\overline{i(X_0)}$ are convex in $\{\uparrow_{O}^{x_0}\}*\Sigma_{\uparrow_{O}^{x_0}}X$.}
Note that there is a geodesic between any
$\uparrow_{O}^{x_1}$ and $\uparrow_{O}^{x_2}$ in $X_0$ (see the beginning of the proof),
then Claim 2 follows from that $i$ is an isometrical embedding.

Due to Claim 2 and ``$X_0\in\mathcal{A}^{n-1}(1)$'', we have $i(X_0),\overline{i(X_0)} \in\mathcal{A}^{n-1}(1)$.

\noindent{\bf Claim 3}: {\it The boundary (in the sense of Alexandrov) $\partial(\overline{i(X_0)})$ is not empty, and $\partial(\overline{i(X_0)})=
\overline{i(X_0)}\setminus i(X_0)$.} Since $\dim(\overline{i(X_0)})=\dim(\{\uparrow_{O}^{x_0}\}*\Sigma_{\uparrow_{O}^{x_0}}X)\ (=n-1)$
and $\overline{i(X_0)}$ is convex in $\{\uparrow_{O}^{x_0}\}*\Sigma_{\uparrow_{O}^{x_0}}X$ (Claim 2),
$\overline{i(X_0)}$ has nonempty boundary in the sense of Alexandrov\footnote{Any
$M\in \mathcal{A}^{n}(k)$ contains no closed convex subset of
dimension $n$ without boundary ([BGP]).}. On the other hand, $X_0$ contains
no boundary point (in the sense of Alexandrov) because $X_0$ is open
in $X$ and $X$ has empty boundary. Since $i$ is an isometrical embedding,
$i(X_0)$ contains no boundary point. It then follows that
$\partial(\overline{i(X_0)})=\overline{i(X_0)}\setminus i(X_0)$. Hence, Claim 3 is verified.

\noindent{\bf Claim 4}: {\it $\overline{i(X_0)}\setminus i(X_0)\subset
\left\{[(\uparrow_{O}^{x_0},w,\frac{\pi}{2})]\in\{\uparrow_{O}^{x_0}\}*\Sigma_{\uparrow_{O}^{x_0}}X\right\}$.}
Note that $i(X_0)$ is open in $\{\uparrow_{O}^{x_0}\}*\Sigma_{\uparrow_{O}^{x_0}}X$ (see the proof of Claim 3), and $|\uparrow_{O}^{x_0}\uparrow_{O}^{x}|<\frac{\pi}{2}$ for any $\uparrow_{O}^{x}\in X_0$ ((i) of Lemma 2.4).
Then it suffices to show that $|\uparrow_{O}^{x_0}v|=\frac{\pi}{2}$ for any
$v\in \overline{X_0}\setminus X_0$ because $i$ is an isometrical embedding.
We first note that $|\uparrow_{O}^{x_0}v|\leq\frac{\pi}{2}$. If $|\uparrow_{O}^{x_0}v|<\frac{\pi}{2}$, then
from the proof of (ii) of Lemma 2.4 it is not hard to see that $\gamma_v$ contains a $x\in N$ with
$|Ox_0|=|Ox|\cos|\uparrow_{O}^{x_0}v|$, and thus $v$
$(=\uparrow_{O}^{x})$ belongs to $X_0$ which contradicts $v\in \overline{X_0}\setminus X_0$. I.e, $|\uparrow_{O}^{x_0}v|=\frac{\pi}{2}$,
so Claim 4 holds.

Due to Claims 3 and 4 and the fact that $\overline{i(X_0)}$ is of dimension
$n-1$, we have
$\overline{i(X_0)}=\{\uparrow_{O}^{x_0}\}*\Sigma_{\uparrow_{O}^{x_0}}X$,
and
$i(X_0)=(\{\uparrow_{O}^{x_0}\}*\Sigma_{\uparrow_{O}^{x_0}}X)^\circ$.
Hence, $X_0$ is isometric to
$(\{\uparrow_{O}^{x_0}\}*\Sigma_{\uparrow_{O}^{x_0}}X)^\circ$ because $i$ is an isometrical embedding.
\hfill $\Box$

\vskip1mm

{\noindent \bf Remark 2.5.1} Consider a cone $C(\{\xi\}*X')$ with vertex $O$, where $X'\in \mathcal{A}(1)$ has empty boundary.
Let $p\in\gamma_{\xi}\setminus\{O\}$. Note that, in the sector $S_{[\xi x]}$ with $x\in X'$, the
geodesic perpendicular to $\gamma_\xi$ at $p$ is parallel to $\gamma_x$ ($\subset \partial(C(\{\xi\}*X'))=C(X')$). Denote by $\beta_x$ such a geodesic, and let
$$N:=\cup_{x\in X'}\beta_x.$$From the proof of Claim 1 in Lemma 2.5, we can see that $N$ is convex in $C(\{\xi\}*X')$,
and $C(\{\xi\}*X')$ is split into two parts $\in \mathcal{A}(0)$ along $N$.

\vskip1mm

{\noindent \bf Definition 2.5.2} In Remark 2.5.1, we say that {\it $N$ is parallel to  $\partial(C(\{\xi\}*X'))$.}

\vskip1mm

From Lemma 2.5 (especially the proof of Claim 1) and Remark 2.5.1, we can draw the following conclusion.

{\noindent \bf Corollary 2.5.3} {\it Suppose that $O\in M_1^\circ$, and
$|Ox_0|_1=\min\{|Ox|_1|x\in N\}$. Then $M_2$ can be isometrically embedded into
$C(\{\uparrow_{O}^{x_0}\}*\Sigma_{\uparrow_{O}^{x_0}}X)$, and $\partial M_2$ is parallel to\\ $\partial(C(\{\uparrow_{O}^{x_0}\}*\Sigma_{\uparrow_{O}^{x_0}}X))$
in $C(\{\uparrow_{O}^{x_0}\}*\Sigma_{\uparrow_{O}^{x_0}}X)$.}

\vskip2mm

{\noindent \bf 2.6\ \ The structure of $X$}

\vskip2mm

Under the condition `$O\in M_1^\circ$', Lemma 2.6 below together with Lemma 2.5 implies that the structure of $X$ has to satisfy one
of the following two cases:

\noindent Case 1: $\overline{X_0}=X$. In this case, $\Sigma_{\uparrow_{O}^{x_0}}X$ admits an isometrical $\Bbb Z_2$-action, which naturally induces
an isometrical $\Bbb Z_2$-action on $S(\Sigma_{\uparrow_{O}^{x_0}}X)$, and $X=S(\Sigma_{\uparrow_{O}^{x_0}}X)/\Bbb Z_2$.

\noindent Case 2: $\overline{X_0}\neq X$. In this case, $\overline{X_0}=\{\uparrow_{O}^{x_0}\}*\Sigma_{\uparrow_{O}^{x_0}}X$ and $X_1:=X\setminus X_0$ is convex in $X$, and thus $X=\overline{X_0}\cup_{\partial\overline{X_0}}X_1$.

\vskip2mm

{\noindent \bf Lemma 2.6} {\it Let $X\in\mathcal{A}^n(1)$ without boundary, and let $X_0$ be open in $X$. Suppose that $X_0$ with the induced metric is isometric to $(\{p\}*X')^\circ$ for some $X'\in\mathcal{A}^{n-1}(1)$ without boundary. Then $X_1:=X\setminus X_0$ is convex in $X$, and one of the following holds:

\noindent{\rm (i)} there is an isometrical $\Bbb Z_2$-action on $X'$, which naturally induces an isometrical $\Bbb Z_2$-action on $S(X')$, and $X=S(X')/\Bbb Z_2$;

\noindent{\rm (ii)} $\overline{X_0}=\{p\}*X'$, and thus $X=\overline{X_0}\cup_{\partial\overline{X_0}}X_1$.}

\vskip1mm

\noindent{\it Proof}. Since $X_0$ with  the induced metric is
isometric to $(\{p\}*X')^\circ$, any $\triangle pq_1q_2$ as a
triangle in $X_0$ is isometric to its comparison triangle in the unit
sphere, and $|px|=\frac{\pi}{2}$ for any
$x\in\overline{X_0}\setminus X_0$ (note that $X_0$ is open in $X$).
It then follows that, for any $x,y\in\overline{X_0}\setminus X_0$,
any geodesic $[xy]$ lies in $X_1$ which implies that $X_1$ is convex
in $X$.

\noindent{\bf Claim 1}: {\it For any $x\in \overline{X_0}\setminus
X_0$, there are at most $2$ geodesics between $p$ and $x$.} In order
to prove the claim, it suffices to show that any two geodesics
$[px]_1$ and $[px]_2$ form an angle equal to $\pi$ at $x$. Assume
that $[px]_1$ and $[px]_2$ form an  angle $\alpha<\pi$ at $x$, which
implies that $|q_1q_2|<|q_1x|+|q_2x|$ for any $q_i\in
[px]_i\setminus\{x\}$ ($i=1$, 2). Note that $|q_ix|\leq |q_iy|$
for any other $y\in X_1$ because $X_0$ is isometric to $(\{p\}*X')^\circ$. Then ``$|q_1q_2|<|q_1x|+|q_2x|$'' implies
that any geodesic $[q_1q_2]$ contains no point in $X_1$ (i.e.
$[q_1q_2]\subset X_0$), and thus $|q_1q_2|=|q_1q_2|_0$, where
$|\cdot|_0$ denotes the induced metric on $X_0$. Now let $q_i$
converge to $x$. Obviously, $|q_1q_2|$ converges to 0. However,
since any $\triangle pq_1q_2$ as a triangle in $X_0$ is isometric to
its comparison triangle, $|q_1q_2|_0$ converges to $\angle q_1pq_2$;
a contradiction. I.e., Claim 1 is verified.

Next we will discuss the structure of $X$ according to the following two cases.

\vskip1mm

\noindent Case 1: $\overline{X_0}=X$ (note that $X_1=\overline{X_0}\setminus X_0$ in this case).

\noindent{\bf Claim 2}: {\it $X_1^2:=\{x\in X_1|\text{ there are $2$ geodesics between $p$ and $x$}\}$
is open and dense in $X_1$.}
We first show that $X_1^2$ is dense in $X_1$. If this is not true, then there exists $x_0\in X_1$ and its neighborhood $U$ in $X_1$ such that there is a unique geodesic between $p$ and any point in $U$. Since $X_1$ is convex in $X$, it follows that the set $V:=\{z\in[px]|x\in U\}$ can be isometrically embedded into $X$ (see Remark A.3.3 in Appendix). Note that $V$ is a neighborhood of $x_0$ in $X$, and thus $x_0$ is a boundary point in the sense of Alexandrov, which contradicts the fact that $X$ has empty boundary. Next we show that $X_1^2$ is open in $X_1$. If this is not true, then there exist $x\in X_1^2$ and a sequence $y_i\in X_1$ such that $y_i\stackrel{i\to\infty}\longrightarrow x$ and there is a unique geodesic between $p$ and $y_i$. Denote by $[px]_j$ ($j=1,2$) the two geodesics between $p$ and $x$. We assume that $[py_i]\longrightarrow[px]_1$ passing a subsequence. By the Subclaim below, for any $q_i\in[px]_2\setminus\{x\}$ and any geodesic $[q_iy_i]$, we have $[q_iy_i]\cap X_1=\{y_i\}$ . This implies that $|q_iy_i|$ converges to the angle between $[px]_1$ and $[px]_2$ at $p$ when $q_i,y_i\longrightarrow x$ (note that $X_0=(\{p\}*X')^\circ$). This obviously yields a contradiction because $|q_iy_i|\longrightarrow 0$.

\noindent{\bf Subclaim}: {\it For any $x\in X_1$, $q\notin X_1$ and any geodesic $[qx]$, we have $[qx]\cap X_1=\{x\}$.} We have proved that there are at most 2 geodesics between $p$ and $x$, and if there are 2 geodesics, then $[px]_1$ and $[px]_2$ form an angle equal to $\pi$ at $x$ (see Claim 1 and its proof). Then by the first variation formula ([BGP]) we obtain that
$$|\uparrow_{x}^p\xi|=\frac\pi2 \text{ for any } \xi\in\Sigma_xX_1 \eqno{(2.8)}$$
(note that  $|px|=\frac{\pi}2$ for all $x\in X_1$ and $X_1$ is
convex in $X$). Since $X_1$ is convex in $X$, if the subclaim is not
true,  then there exists $x'\in X_1$ with $x'\neq x$ such that
$[qx']\setminus\{x'\}\subset X_0$, $[x'x]\subset X_1$ and
$[qx]=[qx']\cup[x'x]$. Note that this implies that $[qx']$ is the
unique geodesic between $q$ and $x'$, and
$$|\uparrow_{x'}^p\uparrow_{x'}^q|+|\uparrow_{x'}^p\uparrow_{x'}^x|=\pi \eqno{(2.9)}$$
for any $\uparrow_{x'}^p$. Note that
$|\uparrow_{x'}^p\uparrow_{x'}^x|=\frac\pi2$  (see (2.8)) because
$[x'x]\subset X_1$. On the other hand, since $X_0=(\{p\}*X')^\circ$
and $[qx']\setminus\{x'\}\subset X_0$, there is a triangle
$\triangle px'q$ which is isometric to its comparison triangle. Then
in this $\triangle px'q$, $\angle px'q<\frac{\pi}{2}$ because
$|px'|=\frac{\pi}{2}$ and $|pq|<\frac{\pi}{2}$, i.e. the geodesic
$[px']$ in this $\triangle px'q$ satisfies
$|\uparrow_{x'}^p\uparrow_{x'}^q|<\frac\pi2$ (note that $[qx']$ is
unique), which contradicts (2.9). Hence, the subclaim is verified.

Since $X_0=(\{p\}* X')^\circ$ (hence $\Sigma_pX=X'$) and $\overline{X_0}=X$, Claim 1 implies that $X'=\{(\uparrow_p^x)_1,(\uparrow_p^x)_2|x\in X_1\}$,
where $(\uparrow_p^x)_i$ is the direction of the geodesic $[px]_i$ at $p$ ($(\uparrow_p^x)_1=(\uparrow_p^x)_2$ if $[px]_1=[px]_2$); and Claim 2 implies that
$X'_2:=\{(\uparrow_p^x)_1,(\uparrow_p^x)_2|x\in X_1^2\}$ is also open and dense in $X'$.  This enables us to define a map
$$\sigma:X'\longrightarrow X'\text{ by } (\uparrow_p^x)_1\longmapsto(\uparrow_p^x)_2 \text{ and } (\uparrow_p^x)_2\longmapsto(\uparrow_p^x)_1 \text{ for any } x\in X_1.$$
Obviously, $\sigma\circ\sigma$ is the identity map. Since $|px|=\frac\pi2$ for any $x\in X_1$ and $X_1$ is convex in $X$, using TCT for `$=$' it is not hard to prove that $\sigma|_{X'_2}$ is a local isometry. Then Sublemma 2.6.1 below and the fact that $X'_2$ is dense in $X'$ imply that $\sigma$ is an isometry. I.e.,
$\sigma$ generates an isometrical $\Bbb Z_2$-action on $X'$.

\noindent{\bf Sublemma 2.6.1}: {\it For any geodesic $[x_0x]$ with
$x_0\in X_1^2$ and $x\in X_1$, $[x_0x]^\circ\subset X_1^2$.
Moreover, for any geodesic $[\uparrow_p^{x_0}\uparrow_p^{x}]\subset
X'$, $[\uparrow_p^{x_0}\uparrow_p^{x}]^\circ$ contains at most one
point which does not belong to $X_2'$.}

\noindent{\it Proof}. Since $|px|=\frac{\pi}2$ for all $x\in X_1$,
due to (2.8) and the TCT for ``='', we can assume that $[px_0]_i,
[x_0x]$ and $[px]_i$ ($i=1,2$ and $[px]_1$ may coincide with
$[px]_2$) form a triangle which is isometric to its comparison
triangle. If $[x_0x]^\circ$ contains a point $y\in X_1\setminus
X_1^2$ (i.e. there is only one geodesic between $p$ and $y$), then
there are two geodesics $[(\uparrow_p^{x_0})_1(\uparrow_p^{x})_1]$,
$[(\uparrow_p^{x_0})_2(\uparrow_p^{x})_2]\subset X'$ which both
contain $\uparrow_p^{y}$. This is impossible because
$$|(\uparrow_p^{x_0})_1(\uparrow_p^{x})_1|=|(\uparrow_p^{x_0})_2(\uparrow_p^{x})_2|=|x_0x|$$
and $$|(\uparrow_p^{x_0})_1(\uparrow_p^{x})_2|,\
|(\uparrow_p^{x_0})_2(\uparrow_p^{x})_1|\geq|x_0x| \text{ (due to
the TCT)}.$$ Moreover, for any $z\in X_1\setminus X_1^2$ (i.e. there
is only one geodesic between $p$ and $z$), the above arguments imply
that any geodesic
$[(\uparrow_p^{x_0})_i\uparrow_p^{z}]=\{\uparrow_p^{x'}|x'\in \text{
some } [x_0z]\}$, and thus
$[(\uparrow_p^{x_0})_i\uparrow_p^{z}]^\circ\subset X_2'$ (because
$[x_0z]^\circ\subset X_1^2$). This implies that for any geodesic
$[\uparrow_p^{x_0}\uparrow_p^{x}]\subset X'$,
$[\uparrow_p^{x_0}\uparrow_p^{x}]^\circ$ contains at most one point
which does not belong to $X_2'$.\hskip1cm $\blacksquare$

Note that $\sigma$ naturally induces an isometry $\bar\sigma$ on
$S(X')$. For convenience, we let $S(X')=\{p,p'\}*X'$ and $\bar
x_i=(\uparrow_p^x)_i$ ($i=1,2$), then we define
$$\bar\sigma: S(X')\longrightarrow S(X')\text{ by } [p\bar x_i]\longmapsto[p'\sigma(\bar x_i)] \text{ and } [p'\bar x_i]\longmapsto[p\sigma(\bar x_i)].$$
Since $\sigma$ is an isometry and $\sigma\circ\sigma=\text{id}$, $\bar\sigma$ is also an isometry and $\bar\sigma\circ\bar\sigma=\text{id}$. And we can define an 1-1 map
$$i:S(X')/\Bbb Z_2\longrightarrow X \text{ by } [[p\bar x_i]]\longmapsto[px]_i,$$
where $[[p\bar x_i]]=\{[p\bar x_i],[p'\sigma(\bar x_i)]\}$.  It is
not hard to check that $i$ is a local isometry, and thus $i$ is an
isometry.

\vskip1mm

\noindent Case 2: $\overline{X_0}\neq X$.

Since $X_1$ is convex in $X$, $X_1\in\mathcal{A}(1)$; and since $X_1\supset X\setminus\overline{X_0}\neq\varnothing$, $X_1$ has the same dimension as $X$.
Thus $X_1$ has nonempty boundary $\partial X_1$ in the sense of Alexandrov geometry (see footnote 6), and obviously $\partial X_1=\overline{X_0}\setminus X_0$.
Note that $|px|=\frac{\pi}{2}$ and $|py|\geq\frac{\pi}{2}$ for any $x\in\partial X_1$ and $y\in X_1$ because $X_0=(\{p\}* X')^\circ$.

\noindent{\bf Claim 3}: {\it For any $x\in \overline{X_0}\setminus X_0$, there is only one geodesic between $p$ and $x$.} Note that there are at most 2 geodesics between $p$ and $x$ (Claim 1). Assume that there are 2 geodesics $[px]_i$ between $p$ and $x$. According to the proof of Claim 1, $[px]_1$ and $[px]_2$ form an angle equal to $\pi$ at $x$, which implies that $\Sigma_xX=\{(\uparrow_x^p)_1,(\uparrow_x^p)_2\}*X''$ for some $X''\in \mathcal{A}^{n-2}(1)$ (see Proposition A.1 in Appendix and note that $\dim(\Sigma_xX)=n-1$).
On the other hand, since $|px|=\frac{\pi}{2}$ and $|py|\geq\frac{\pi}{2}$ for any $x\in\partial X_1$ and $y\in X_1$ and $X_1$ is convex in $X$, by the first variation formula, $$|(\uparrow_x^p)_i\xi|\geq\frac{\pi}{2} \text{ for any } \xi\in \Sigma_x X_1.\eqno{(2.10)}$$
Since $\Sigma_xX=\{(\uparrow_x^p)_1,(\uparrow_x^p)_2\}*X''$, (2.10) implies that $\Sigma_x X_1\subseteq X''$. However, $\dim(\Sigma_x X_1)=\dim(\Sigma_x X)=n-1$ and $\dim(X'')=n-2$,
a contradiction; i.e., Claim 3 is true.

\noindent{\bf Claim 4}: {\it $\overline{X_0}$ and $\overline{X_0}\setminus X_0$ with the induced metrics from $X$ are isometric to $\{p\}* X'$ and $X'$ respectively.}
Note that ``$X_0=(\{p\}* X')^\circ$'' implies that the shortest path in $\overline{X_0}$ between any two points in $\overline{X_0}\setminus X_0$ still falls in $\overline{X_0}\setminus X_0$. Then by the first variation formula ([BGP]), ``$|px|=\frac{\pi}{2}$ for all $x\in \overline{X_0}\setminus X_0=\partial X_1$'' and Claim 3 imply that $|\uparrow_x^p\xi|=\frac{\pi}{2}$ for any $\xi\in\Sigma_x(\partial X_1)$ (note that $X_1\in\mathcal{A}(1)$). Similar to getting the subclaim in the proof of Claim 2, we can conclude that $[qx]_{\overline{X_0}}\cap\partial X_1=\{x\}$ for any $q\in X_0$ and $x\in\partial X_1$, and thus the shortest path in $\overline{X_0}$ between any two points in $X_0$ still falls in $X_0$. Then Claim 4 follows from ``$X_0=(\{p\}* X')^\circ$''.

On the other hand, since $X_1$ is convex in $X$, the induced metric on $\partial X_1$ from $X_1$ is just that from $X$. This together with Claim 4 implies that
$\partial\overline{X_0}$ is isometric to $\partial X_1$ with respect to the induced metrics from $\overline{X_0}$ and $X_1$ respectively. Hence, we conclude that $X=\overline{X_0}\cup_{\partial\overline{X_0}}X_1$ because $\overline{X_0}=\{p\}* X'$ and $X_1$ is convex in $X$.
\hfill $\Box$

\vskip2mm

{\noindent \bf 2.7\ \ On the case ``$O\in N$''}

\vskip2mm

{\noindent \bf Lemma 2.7} {\it Suppose that $O\in N$. Then there exist $X_i\in \mathcal{A}^{n-1}(1)$ with boundaries $\partial X_1\stackrel{\rm{iso}}{\cong}\partial X_2$ such that $X=X_1\cup_{\partial X_i}X_2$ and $M_i=C(X_i)$.}

\vskip1mm

\noindent{\it Proof}. We firstly observe that for any $p\in M_i^\circ$, the geodesic $[Op]_{M_i}$ is just $[Op]_M$. In fact, $[Op]\setminus\{O\}$ is a local geodesic in $M$ because $[Op]\setminus\{O\}\subset M_i^\circ$ (see footnote 5), so by Proposition 2.7.1 below $[Op]_{M_i}$ has to lie in the ray $\gamma_{\uparrow_O^p}$ in $M$, i.e. $[Op]_{M_i}=[Op]_M$.

{\noindent \bf Proposition 2.7.1} {\it Let $C(X)$ be a cone with vertex $O$ for some $X\in \mathcal{A}(1)$, and let $c$ be a local geodesic on $C(X)$.
Then either $c$ lies in some ray starting from $O$ or there is $\delta>0$ such that $|Op|\ge\delta$ for any $p\in c$.}

\noindent{\it Proof}. If $c$ does not lie in any ray starting from $O$, then $c$ is a geodesic in the (Euclidean) sector $S_c$ ($=\{\gamma_{\uparrow_O^p}|p\in c\}$). Moreover,
there is $\delta>0$ such that the distance between $O$ and $c$ is not less than $\delta$ in the sector $S_c$. Note that $|Oq|_{S_c}=|Oq|_{C(X)}$ for any $q\in S_c$. Then it follows that $|Op|\ge\delta$ for any $p\in c$.  \hskip1cm $\blacksquare$

Due to the above observation, we claim that any geodesic $\gamma$ in $M$ starting from $O$ lies in $M_1$ or $M_2$. In fact, if the claim is not true, then there are $p_1,p_2\in \gamma$
with $p_i\in M_i^\circ$. It then follows from the above observation that $[Op_i]_{M_i}=[Op_i]_M\subset\gamma$ for $i=1$ and $2$. Without loss of generality, we assume that $|Op_1|<|Op_2|$. Then $p_1\in [Op_2]^\circ_M=[Op_2]^\circ_{M_2}\subset M_2^\circ$ which contradicts $p_1\in M_1^\circ$.

Due to the claim, $\Sigma_OM_i=\{(\uparrow_O^x)_M|x\in M_i\}$ and $M_i=C(\Sigma_OM_i)$ for $i=1$ and 2. Of course, $\Sigma_OM_i$ has nonempty boundary  because $M_i$ has nonempty boundary. From the definition of the metric of the cone, it is not hard to see that the metric of $\Sigma_OM_i$ is just the induced metric from $X=\Sigma_OM$, and $\partial(\Sigma_OM_1)\stackrel{\rm{iso}}{\cong}\partial(\Sigma_OM_2)$ because $\partial M_1\stackrel{\rm{iso}}{\cong}\partial M_2$.
It then follows that $$X=\Sigma_OM_1\cup_{\partial(\Sigma_OM_i)}\Sigma_OM_2.$$
Hence the proof is done once we let $X_i=\Sigma_OM_i$  for $i=1$ and 2.
\hfill $\Box$

%%%%%%%%%%%%%%%%%%%%%%%%%%%%%%%%%%%% Section 3 %%%%%%%%%%%%%%%%%%%%%%%%%%%%%%%%%%%%%%%%

\section{Proof of Theorem C}

We will apply induction on $\dim(Y_1*Y_2)$ to prove Theorem C. When $\dim(Y_1)=0$, ``$Y_1$ has empty boundary'' means that $Y_1=\{p,q\}$ with $|pq|=\pi$, and thus $Y_1*Y_2=S(Y_2)$.
Hence, if $\dim(Y_1)=0$, then Theorem C follows from Corollary 0.1. In the rest of this section, we assume that $\dim(Y_1)>0$ and $\dim(Y_2)>0$.

We will give the proof according to two cases: the diameter $\diam(Y_1*Y_2)=\pi$ and $\diam(Y_1*Y_2)<\pi$. Subsection 3.1 is on the former case, and Subsection 3.6 is on the latter case
(Subsections 3.2-3.5 are prepared for 3.6).

\vskip2mm

{\noindent \bf 3.1\ \ On the case that the diameter of $Y_1*Y_2$ is $\pi$}

\vskip2mm

{\noindent \bf Lemma 3.1} {\it Theorem C holds if the diameter $\diam(Y_1*Y_2)=\pi$.}

\vskip2mm

\noindent{\it Proof}. Since $\diam(Y_1*Y_2)=\pi$, there are $p_1,p_2\in Y_1*Y_2$ with $|p_1p_2|=\pi$ and $Z\in \mathcal{A}(1)$ such that
$Y_1*Y_2=\{p_1,p_2\}*Z=S(Z)$ (see Proposition A.1). On the other hand, $\diam(Y_1*Y_2)=\pi$ implies that $\diam(Y_1)=\pi$ because $\diam(Y_2)<\pi$ (Proposition A.5),
so $p_1,p_2\in Y_1$ and there exists $Y_1'\in \mathcal{A}(1)$ such that $Y_1=\{p_1,p_2\}*Y_1'$ and $Z=Y_1'*Y_2$.

By Corollary 0.1, $p_i$ and $Z$ can be chosen such that $M_i=\{p_i\}*Z$ or $M_i=\{p_1,p_2\}*Z_i$ where $Z_1\cup_{\partial Z_i}Z_2=Z$, which correspond to the following two cases respectively:

\noindent Case 1: $M_i=(\{p_i\}*Y_1')*Y_2$, and we only need to let $X_i=\{p_i\}*Y_1'$.

\noindent Case 2: $Z_1\cup_{\partial Z_i}Z_2=Z=Y_1'*Y_2$. By induction, there exist $X_i'\in \mathcal{A}(1)$ with $\partial X'_1\stackrel{\rm{iso}}{\cong}\partial X'_2$ such that $Y_1'=X'_1\cup_{\partial X'_i}X'_2$ and $Z_i=X'_i*Y_2$, or $Y_2=X'_1\cup_{\partial X'_i}X'_2$ and $Z_i=Y_1'*X_i'$; and
we only need to let $X_i=\{p_1,p_2\}* X'_i$ or $X_i'$ respectively.
\hfill $\Box$

\vskip2mm

In the rest of the proof, we only need to solve the case ``$\diam(Y_1*Y_2)<\pi$'' (i.e. each $\diam(Y_j)<\pi$). In the following, Subsections 3.2 and 3.4-3.6 are under the conditions of Theorem C and $\diam(Y_1*Y_2)<\pi$; and we still use $M$ and $N$ to denote $Y_1*Y_2$ and $\partial M_i$ respectively; and $Y_1$ and $Y_2$ as subsets of $Y_1*Y_2$ are sets $\{[(y_1,*,0)]\in Y_1*Y_2|y_1\in Y_1\}$ and $\{[(*,y_2,\frac{\pi}{2})]\in Y_1*Y_2|y_2\in Y_2\}$ respectively.

\vskip2mm

{\noindent \bf 3.2\ \ The direction space at a point in $N\cap Y_i$}

\vskip2mm

Let $p$ be a point in $N\cap Y_1$ or $N\cap Y_2$, say $N\cap Y_1$. It then follows from the gluing theorem that $\Sigma_pM=\Sigma_pM_1\cup_{\partial(\Sigma_pM_i)}\Sigma_pM_2$ (ref. [Pet]). On the other hand, $\Sigma_pM=Y_1'*Y_2'$ where $Y_1'=\Sigma_pY_1$ and $Y_2'=Y_2$ (see Proposition A.4). Note that $\diam(Y_2')<\pi$, then by the induction
there are $X_1',X_2'\in \mathcal{A}(1)$ with $\partial X_1'\stackrel{\rm{iso}}{\cong}\partial X_2'$ such that $Y_1'=X'_1\cup_{\partial X'_i}X'_2$ and $\Sigma_pM_i=X'_i*Y'_2$, or $Y'_2=X'_1\cup_{\partial X'_i}X'_2$ and $\Sigma_pM_i=Y_1'*X_i'$; and $\partial(\Sigma_pM_i)=\partial X_i'*Y'_2$ or $Y'_1*\partial X_i'$ respectively.

\vskip2mm

{\noindent \bf 3.3\ \ A basic proposition of the join $Y_1*Y_2$}

\vskip2mm

From the definition of the metric of the join, for any geodesics
$[y_jy_j']\subset Y_j$ $(j=1,2)$, $[y_1y_1']*[y_2y_2']$ can be
isometrically embedded into $Y_1*Y_2$ (also into the unit 3-sphere $\Bbb
S^3=S^1*S^1$ where the two $S^1$ have diameter equal to $\pi$). Of
course, for a geodesic $\gamma:[0,\ell]\to
[y_1y_1']*[y_2y_2']\subset Y_1*Y_2$ defined by $s\mapsto
[(y_1(s),y_2(s),t(s))]$ with $\gamma(0)\in[y_1y_2]^\circ$ and
$\gamma(\ell)\in[y_1'y_2']^\circ$, $y_j(s)|_{[0,\ell]}$ is just the
geodesic $[y_jy_j']$ for $j=1$ and 2. How about any geodesic
$\gamma:[0,\ell]\to\subset Y_1*Y_2$?

\vskip2mm

{\noindent \bf Proposition 3.3}\ \ {\it Let $Y_1,Y_2\in
\mathcal{A}(1)$, and let $\gamma:[0,\ell]\to Y_1*Y_2$ defined by
$s\mapsto [(y_1(s),y_2(s),t(s))]$ be a geodesic with
$t(0), t(\ell)\neq0$ or $\frac\pi2$. Then $y_j(s)|_{[0,\ell]}$ is a
geodesic in $Y_j$ for $j=1$ and $2$, and $\gamma\subset
y_1(s)|_{[0,\ell]}*y_2(s)|_{[0,\ell]}\subset Y_1*Y_2$.}

\vskip2mm

\noindent{\it Proof}. We consider $\gamma|_{[0,s_0]}$ with $s_0<\ell$. Note that $\gamma|_{[0,s_0]}$ is the unique
geodesic between $\gamma(0)$ and $\gamma(s_0)$. On the other hand, from the definition of the metric of the join,
there is a geodesic between $\gamma(0)$ and $\gamma(s_0)$ in $[y_1(0)y_1(s_0)]*[y_2(0)y_2(s_0)]\subset Y_1*Y_2$ for any
any geodesic  $[y_1(0)y_1(s_0)]$ and $[y_2(0)y_2(s_0)]$. Hence, there is a unique geodesic between $y_j(0)$ and $y_j(s_0)$.
Due to the arbitrary of $s_0$ and the continuity, $y_j(s)|_{[0,\ell]}$ is a geodesic in $Y_j$, and $\gamma\subset
y_1(s)|_{[0,\ell]}*y_2(s)|_{[0,\ell]}\subset Y_1*Y_2$.
\hfill $\Box$

\vskip2mm

{\noindent \bf Remark 3.3.1}\ \ If $\gamma$ in Proposition 3.3 satisfies $t(0)=0$ and
$t(\ell)\neq0$ or $\frac\pi2$, then $y_1(s)|_{[0,\ell]}$ is a
geodesic in $Y_1$, and $\gamma\subset
y_1(s)|_{[0,\ell]}*\{y_2(\ell)\}\subset Y_1*Y_2$.

\vskip2mm

{\noindent \bf Remark 3.3.2}\ \ {In Proposition 3.3, if $y_j(0)=y_j(\ell)$ for $j=1$ or $2$, then the conclusion of the proposition implies that $y_j(s)|_{[0,\ell]}\equiv y_j(0)=y_j(\ell)$.}

\vskip2mm

From Proposition 3.3 and Remark 3.3.1, it is not hard to see the following corollary.

\vskip2mm

{\noindent \bf Corollary 3.3.3}\ \ {\it Let $Y_1,Y_2\in
\mathcal{A}(1)$, and let $\gamma:(0,\ell)\to Y_1*Y_2$  be a local geodesic parameterized by arc length $s$ and defined by
$s\mapsto [(y_1(s),y_2(s),t(s))]$ with $t(s)\not\equiv0$ or $\frac\pi2$.
Then  $y_j(s)|_{(0,\ell)}$ is a local
geodesic in $Y_j$ for $j=1$ and $2$; and there is a local geodesic $\tilde\gamma:(0,\ell)\to \Bbb S^3=S^1* S^1$  also parameterized by arc length $s$ and defined by
$s\mapsto [(\theta_1(s),\theta_2(s),\tilde t(s))]$ such that the local isometry
$$i:\gamma\to\tilde\gamma \text{ defined by } \gamma(s)\mapsto \tilde\gamma(s)$$ satisfies that $y_j(s)|_{(0,\ell)}\to \theta_j(s)|_{(0,\ell)}$ is a local
isometry and $\tilde t(s)=t(s)$ for any $s$.}

\vskip2mm

{\noindent \bf Remark 3.3.4}\ \ {In Corollary 3.3.3, $t(s)$ has the following three cases:

\noindent Case 1: There is $\delta>0$ such that $t(s)\in[\delta, \frac\pi2-\delta]$ for any $s\in(0,\ell)$.

\noindent Case 2: There is $s_0$ such that $t(s_0)=0$ (resp. $t(s_0)=\frac{\pi}{2}$), and $t(s)<\frac\pi2$ (resp. $t(s)>0$)  for any $s$. In this case, $\{s\in (0,\ell)|t(s)=0 \text{ (resp. $t(s)=\frac{\pi}{2}$) }\}=\{s_1, s_2, \cdots, s_k\}$ with
$s_{i+1}-s_{i}=\pi$ (of course $k$ may be equal to 1), and $y_2(s)|_{(0,s_1)}$, $y_2(s)|_{(s_i,s_{i+1})}$ and $y_2(s)|_{(s_{k},\ell)}$ (resp. $y_1(s)|_{(0,s_1)}$,  $y_1(s)|_{(s_i,s_{i+1})}$ and $y_1(s)|_{(s_{k},\ell)}$) are isolated points respectively (we also say that $y_j(s)|_{(0,\ell)}$ is a local geodesic).

\noindent Case 3: $A:=\{s\in (0,\ell)|t(s)=0 \text{ or } \frac{\pi}{2}\}=\{s_1, t_1, s_2, t_2,\cdots, s_k, t_k|t(s_i)=0, t(t_i)=\frac{\pi}{2})\}$ with
$s_{i+1}-t_{i}=t_i-s_i=\frac{\pi}{2}$ ($A$ maybe have the forms $\{s_1,t_1, \cdots, s_k\}$, or $\{t_1, s_2, t_2,\cdots, s_k, t_k\}$, or $\{t_1, s_2, t_2,\cdots, s_k\}$). In this case, $\gamma|_{(0,s_1]}\subseteq [y_2(0)y_1(s_1)]$, $\gamma|_{[s_i,t_i]}=[y_1(s_i)y_2(t_i)]$, $\gamma|_{[t_i,s_{i+1}]}=[y_2(t_i)y_1(s_{i+1})]$ and $\gamma|_{[t_k,\ell)}\subset[y_2(t_k)y_1(\ell)]$. }

\vskip2mm

For convenience of readers, we give a detailed proof for an easy case of Corollary 3.3.3.

{\noindent \bf An easy case of Corollary 3.3.3}:\ \ {\it In Corollary 3.3.3, if $s\in[0,\ell]$ and if $t(0)=0$ (i.e. $\gamma(0)\in Y_1\subset
Y_1*Y_2$) and $0<t(s)<\frac{\pi}{2}$ for any $s\in(0,\ell]$, then $y_1(s)|_{[0,\ell]}$ is a local geodesic in $Y_1$,
and $\gamma$ is a geodesic in
$y_1(s)|_{[0,\ell]}*\{y_2(s_0)\}\subset Y_1*Y_2$ for some
$s_0\in(0,\ell)$.}

\vskip1mm

Note that $y_1(s)|_{[0,\ell]}*\{y_2(s_0)\}$ can be isometrically embedded into $\Bbb S^3=S^1*S^1$.

\vskip1mm

\noindent{\it Proof}. Since $\gamma$ is a local geodesic, there
exists $s_0>0$ such that $\gamma|_{[0,s_0]}$ is a geodesic. By Remark 3.3.1, $y_1(s)|_{[0,s_0]}$ is a
geodesic in $Y_1$, and $\gamma|_{[0,s_0]}$ is a geodesic in $y_1(s)|_{[0,s_0]}*\{y_2(s_0)\}\subset Y_1*Y_2$.
Furthermore, we can select $s_1,s_2\in[0,\ell]$ such that
$s_1<s_0<s_2$ and $\gamma|_{[s_1,s_2]}$ is the unique geodesic
between $\gamma(s_1)$ and $\gamma(s_2)$.

\noindent{\bf Claim}: {\it $y_1(s)|_{[s_1,s_2]}$ is a geodesic in $Y_1$, and
$\gamma|_{[0,s_2]}$ is a geodesic in $\{y_2(s_0)\}*y_1(s)|_{[0,s_2]}$.}

By assuming that the claim is true, we can draw the conclusion step by step (next step is to select
$s_3,s_4\in[0,\ell]$ with $s_3<s_2<s_4$ such that
$y_1(s)|_{[s_3,s_4]}$ is a geodesic in $Y_1$ and $\gamma|_{[0,s_4]}$
is a geodesic in $y_1(s)|_{[0,s_4]}*\{y_2(s_0)\}$).

In the rest of the proof, we only need to verify the claim. Since
$\gamma|_{[s_1,s_2]}$ is the unique geodesic between $\gamma(s_1)$
and $\gamma(s_2)$, triangles
$\triangle y_2(s_0)\gamma(s_1)\gamma(s_2)$, $\triangle
y_2(s_0)\gamma(s_1)\gamma(s_0)$ and $\triangle
y_2(s_0)\gamma(s_0)\gamma(s_2)$ in $Y_1*Y_2$ are isometrical to
their comparison triangles respectively (see Corollary A.3.1 in Appendix), which implies that
$$|\uparrow_{y_2(s_0)}^{\gamma(s_1)}\uparrow_{y_2(s_0)}^{\gamma(s_2)}|
=|\uparrow_{y_2(s_0)}^{\gamma(s_1)}\uparrow_{y_2(s_0)}^{\gamma(s_0)}|+
|\uparrow_{y_2(s_0)}^{\gamma(s_0)}\uparrow_{y_2(s_0)}^{\gamma(s_2)}|
=|y_1(s_1)y_1(s_0)|+|\uparrow_{y_2(s_0)}^{\gamma(s_0)}\uparrow_{y_2(s_0)}^{\gamma(s_2)}|
\eqno{(3.1)}$$ and thus
$$|\uparrow_{y_2(s_0)}^{\gamma(s_1)}\uparrow_{y_2(s_0)}^{\gamma(s_2)}|
-|\uparrow_{y_2(s_0)}^{\gamma(s_0)}\uparrow_{y_2(s_0)}^{\gamma(s_2)}|
\geq|y_1(s_1)y_1(s_2)|-|y_1(s_0)y_1(s_2)|;\eqno{(3.2)}$$ and by Remark A.4.2 in Appendix
$$\begin{aligned}& \cos|\uparrow_{y_2(s_0)}^{\gamma(s_0)}\uparrow_{y_2(s_0)}^{\gamma(s_2)}|=\cos|y_1(s_0)y_1(s_2)|\cos|\uparrow_{y_2(s_0)}^{y_1(s_2)}\uparrow_{y_2(s_0)}^{\gamma(s_2)}|,\\
&\cos|\uparrow_{y_2(s_0)}^{\gamma(s_1)}\uparrow_{y_2(s_0)}^{\gamma(s_2)}|=\cos|y_1(s_1)y_1(s_2)|\cos|\uparrow_{y_2(s_0)}^{y_1(s_2)}\uparrow_{y_2(s_0)}^{\gamma(s_2)}|.
\end{aligned}\eqno{(3.3)}$$

\noindent{\bf Sublemma 3.3.5}: {\it For any angles $\alpha_i,
\beta_i, \theta\in[0,\frac\pi2]$ $(i=1,2)$ with
$0\leq\beta_2-\beta_1\leq\alpha_2-\alpha_1$, if
$\cos\alpha_i=\cos\beta_i\cos\theta$ for $i=1$ and $2$, then
$\alpha_i=\beta_i$ and $\theta=0$.}

\noindent{\it Proof}. The proof is omitted because it is
elementary.\hskip1cm $\blacksquare$

By Sublemma 3.3.5, if in addition $s_1$ and $s_2$ are selected such
that $|\uparrow_{y_2(s_0)}^{\gamma(s_1)}\uparrow_{y_2(s_0)}^{\gamma(s_2)}|\leq\frac\pi2$ and $|y_1(s_0)y_1(s_2)|\leq|y_1(s_1)y_1(s_2)|\leq\frac\pi2$,
then (3.2) and (3.3) together implies that
$$|\uparrow_{y_2(s_0)}^{\gamma(s_1)}\uparrow_{y_2(s_0)}^{\gamma(s_2)}|
=|y_1(s_1)y_1(s_2)|,\
|\uparrow_{y_2(s_0)}^{\gamma(s_0)}\uparrow_{y_2(s_0)}^{\gamma(s_2)}|=|y_1(s_0)y_1(s_2)|,\eqno{(3.4)}$$
and
$$|\uparrow_{y_2(s_0)}^{y_1(s_2)}\uparrow_{y_2(s_0)}^{\gamma(s_2)}|=0.\eqno{(3.5)}$$
Obviously, (3.1) and (3.4) imply that $$|y_1(s_1)y_1(s_2)|=|y_1(s_1)y_1(s_0)|+|y_1(s_0)y_1(s_2)|.\eqno{(3.6)}$$
It is easy to check that the above process implies that $|y_1(s_1)y_1(s_2')|=|y_1(s_1)y_1(s_0)|+|y_1(s_0)y_1(s_2')|$
for any $s_2'\in(s_0,s_2)$, i.e. (3.6) implies that $y_1(s)|_{[s_1,s_2]}$ is a
geodesic in $Y_1$. On the other hand, (3.5) implies that
$\gamma(s_2)\in[y_2(s_0)y_1(s_2)]$ (i.e $y_2(s_2)=y_2(s_0)$), and
thus $y_2(s)|_{[s_1,s_2]}=y_2(s_0)$ and $\gamma|_{[0,s_2]}$ is a
geodesic in $y_1(s)|_{[0,s_2]}*\{y_2(s_0)\}$. \hfill $\Box$

%%%%%%%%%%%%%%%%%%%%%%%%%%%%%%%%%%%%%%%%%%%%%%%%%%%%%%%%%%%%%%%%%%%%%%%%%%%%%%%%%%%%%%%%%%%

\vskip2mm

{\noindent \bf 3.4\ \ About $[p_1p_2]_M\cap N$ for any $p_j\in
Y_j\subset Y_1*Y_2$}

\vskip2mm

{\noindent \bf Lemma 3.4} {\it For any $p_j\in Y_j\subset Y_1*Y_2$ with  $j=1$ and $2$,

\noindent{\rm (i)} if $p_1,p_2\in M_i$ and $p_1\in M_i^\circ$, then $[p_1p_2]_M\subset M_i$;

\noindent{\rm (ii)} if $p_1,p_2\in N$, then $[p_1p_2]_M\subset N$;

\noindent{\rm (iii)} if $p_1,p_2\not\in N$, then  $[p_1p_2]_M$ contains at most one point in $N$.}

\vskip2mm

\noindent{\it Proof}. (i)\ \ Since $p_j\in M_i$ and $p_1\in M_i^\circ$, $[p_1p_2]_{M_i}\setminus\{p_2\}$ is a local geodesic in $M$.
From Corollary 3.3.3 (see Case 3 in Remark 3.3.4), the length $\ell([p_1p_2]_{M_i})=\frac\pi2$ or $\geq\frac{3\pi}{2}$. We know that
$\ell([p_1p_2]_{M_i})\leq\pi$ because $[p_1p_2]_{M_i}$ is a geodesic in $M_i\in\mathcal{A}(1)$.
Hence, $\ell([p_1p_2]_{M_i})=\frac\pi2$, and thus $[p_1p_2]_{M_i}$ has to be $[p_1p_2]_M$.

\vskip1mm

(ii)\ \ Note that in this case, only the following two subcases may occur.

Subcase 1: For any neighborhoods $U_j\subset Y_j$ of $p_j$, $U_1\cup U_2$  contains points in $M_1^\circ$ and $M_2^\circ$.
By (i) and the limit argument, we conclude that $[p_1p_2]_{M}\subset M_1\cap M_2=N$.

Subcase 2: For $j=1$ and 2, there are neighborhoods $U_j\subset Y_j$ of $p_j$ such that both $U_1$ and $U_2$ belong to $M_1$ or $M_2$, say $M_1$.

In this case, we first prove that $[p_1p_2]_{M_1}\subset N$.
In fact, if $[p_1p_2]_{M_1}\not\subset N$, then $[p_1p_2]^\circ_{M_1}\subset M_1^\circ$ (see footnote 5). This implies that
there are neighborhoods $V_j\subset U_j$ of $p_j$ such that $[p_1'p_2']^\circ_{M_1}\subset M_1^\circ$ for any $p'_j\in V_j$ (note that
either $[p_1'p_2']^\circ_{M_1}\subset M_1^\circ$ or $[p_1'p_2']^\circ_{M_1}\subset N$), and thus
$[p_1'p_2']^\circ_{M_1}$ is a local geodesic in $M$.
By the same argument as proving (i), we conclude that $[p_1'p_2']_{M_1}=[p_1'p_2']_M$.
I.e., $V_1*V_2\subset M_1$ and $(V_1*V_2)\setminus (V_1\cup V_2)\subset M_1^\circ$, which implies that
$\dim(\partial(\Sigma_{p_j}M_1))\leq \dim(Y_j)$. However, $\dim(\partial(\Sigma_{p_j}M_1))=\dim(M)-1=\dim(Y_1)+\dim(Y_2)>\dim(Y_j)$
(note that we have assumed $\dim(Y_j)>0$), which contradicts $\dim(\partial(\Sigma_{p_j}M_1))\leq \dim(Y_j)$.

Since $[p_1p_2]_{M_1}\subset N$, we only need to show that $[p_1p_2]_{M_1}=[p_1p_2]_M$. Note that only the following two cases may occur:
$W_1\cup W_2$ contains points in $M_1^\circ$ for any neighborhoods $W_j\subset U_j$ of $p_j$;
there are neighborhoods $W_j\subset U_j$ of $p_j$ such that $W_j\subset N$ ($j=1,2$).
In the former case, we conclude that $[p_1p_2]_{M}\subset M_1$ like in Subcase 1, and thus $[p_1p_2]_{M_1}=[p_1p_2]_M$.
In the latter case, we will derive a contradiction (and the proof is done). We first give a claim.
{\bf Claim 1}: {\it $Y_j\subset N$ for $j=1$ and $2$.} If $Y_j$ contains a point $y_j$ belonging to $M_i^\circ$,
then $[p_jy_j]_{M_i}\setminus\{p_j\}\subset M_i^\circ$ (see footnote 5) which implies that $[p_jy_j]_{M_i}\setminus\{p_j\}$ is a local geodesic in $M$.
Moreover, since $W_j\subset N$, $[p_jy_j]_{M_i}\not\subset Y_j$. By Corollary 3.3.3 (and Remark 3.3.4), the length of $[p_jy_j]_{M_i}\setminus\{p_j\}$ is equal to $\pi$, i.e. $|p_jy_j|_i=\pi$. Hence, such $y_j$ is unique, i.e. $Y_j\cap M_i^\circ=\{y_j\}$.
However, ``$y_j\in M_i^\circ$'' implies that there is a neighborhood $U\subset Y_j$ of $y_j$ such that $U\subset M_i^\circ$; a contradiction.  I.e., Claim 1 is verified.
Note that Claim 1 implies that $Y_1$ and $Y_2$ are convex in $M_i$. Since $Y_j$ has empty boundary, by Sublemma 3.4.1 below $|y_1y_2|_i=\frac\pi2$ for any
$y_j\in Y_j$ (note that $|y_1y_2|_i\geq |y_1y_2|=\frac\pi2$). This implies that $[y_1y_2]_M\subset M_i$, i.e., $M_i=M$; a contradiction.

\vskip1mm

\noindent {\bf Sublemma 3.4.1 ([Ya]).} {\it Let $X\in\mathcal{A}(1)$, and let $Y$ be a complete locally convex subset
without boundary in $X$. If $|py|\geq\frac{\pi}{2}$ for some $p\in X$ and any $y\in Y$, then $|py|=\frac{\pi}{2}$.}

\vskip1mm

(iii)\ \ Assume that $p_1\in M_1^\circ$. We consider the following
set $$X=\left\{[(y_1, y_2, t)]\in M|y_1\in Y_1\cap M_1 \text{ and }
[(y_1,y_2,t')]|_{0\leq t'\leq t}\subset M_1 \right\}.$$ Obviously,
$X$ is a closed subset in $M$. {\bf Claim 2}: {\it $\partial X$ $(=X\setminus X^\circ)$
belongs to $N$.} If the claim is not true, then there exists
$x=[(y_1,y_2,t_0)]\in\partial X$ such that $x\in M_1^\circ$. If
$y_1\in M_1^\circ$, then it is not hard to see that there exists
$B_x(\epsilon)\subset M_1^\circ$ belonging to $X$ (note that
$[y_1x]_M\subset M_1^\circ$), i.e. $x$ is an inner point of $X$; a
contradiction. Hence, it has to hold that $y_1\in N$. Since $p_1\in
M_1^\circ$, $[y_1p_1]_{M_1}\setminus\{y_1\}\subset M_1^\circ$ (see
footnote 5) and thus $[y_1p_1]_{M_1}\setminus\{y_1\}$ is a local
geodesic in $M$. By Corollary 3.3.3 (and Remark 3.3.4),
$[y_1p_1]_{M_1}$ has to belong to $Y_1$ (otherwise the length
$\ell([y_1p_1]_{M_1})=\pi$ (i.e. $|y_1p_1|_1=\pi$), and thus
$M_1=\{p_1,y_1\}*A$ for some $A\in \mathcal{A}(1)$ with $\partial
M_1=\{p_1,y_1\}*\partial A$ (see Corollary A.4.1 in Appendix) which contradicts ``$p_1\in
M_1^\circ$''). Hence,
$(\uparrow_{y_1}^{p_1})_{M_1}\in(\Sigma_{y_1}M_1)^\circ\cap\Sigma_{y_1}Y_1$;
and then by the induction on $\Sigma_{y_1}M=(\Sigma_{y_1}Y_1)*Y_2$
(see 3.2), there are $X_1',X_2'\in \mathcal{A}(1)$ such that
$\Sigma_{y_1}Y_1=X'_1\cup_{\partial X'_i}X'_2$ and
$\Sigma_pM_i=X'_i*Y_2$ and $\partial(\Sigma_pM_i)=\partial
X_i'*Y_2$. It follows that
$(\uparrow_{y_1}^{y_2'})_M\in\partial(\Sigma_{y_1}M_i)$ for any
$y_2'\in Y_2$, which contradicts
``$(\uparrow_{y_1}^{y_2})_M=(\uparrow_{y_1}^{x})_{M_1}\in(\Sigma_{y_1}M_1)^\circ$''.
That is, Claim 2 is verified.

Now we consider $[p_1p_2]_M\cap N$ and assume that $[p_1p_2]_M\cap
N\neq\emptyset$. Let $z$ be the first point in $N$ along
$[p_1p_2]_M$ (from $p_1$ to $p_2$). It suffices to show that
$[zp_2]_M\cap N=\{z\}$. Assume that there is $z'\in [zp_2]_M\cap N$
with $z'\neq z$. {\bf Claim 3}: $z'\not\in X$. If $z'\in X$, then
$[p_1z']_M\subset M_1$ which implies that $[p_1z']_M=[p_1z']_{M_1}$.
Since $p_1\in M_1^\circ$, $[p_1z']_M\setminus\{z'\}\subset
M_1^\circ$ (see footnote 5), which contradicts ``$z\in N$'' (note
that $z\in[p_1z']_M\setminus\{z'\}$). On the

\noindent other hand, $p_1\in X^\circ$ because $p_1\in M_1^\circ$
(see Claim 2). It then follows that $[p_1z']^\circ_{M_1}\cap\partial
X\neq\emptyset$. This is impossible because
$[p_1z']^\circ_{M_1}\subset M_1^\circ$ and $\partial X\subset N$
(see Claim 2). Hence, it has to hold that $[zp_2]_M\cap N=\{z\}$, and the proof is done.
\hfill $\Box$

\vskip2mm

{\noindent \bf 3.5\ \ The relations between $N$ and $Y_j$}

\vskip2mm

{\noindent \bf Lemma 3.5} {\it $N\cap Y_j\neq\emptyset$ for $j=1$
and $2$. Moreover, either $Y_1\cap M_i^\circ\neq\emptyset$ ($i=1,
2$) and $Y_2\subset N$, or $Y_2\cap M_i^\circ\neq\emptyset$ ($i=1,
2$) and $Y_1\subset N$.}

\vskip2mm

{\noindent \bf Sublemma 3.5.1} {\it If $y_1\in Y_1\cap M_1^\circ$,
then $(\{y_1\}*Y_2)\cap N\subset Y_2 (\subset Y_1*Y_2)$ which
implies that $\{y_1\}*Y_2\subset M_1$, where $\{y_1\}*Y_2\subset
Y_1*Y_2$.}

\vskip1mm

\noindent{\it Proof}. We give the proof according to the following
three cases.

Case 1: $Y_2\subset M_1$. In this case, by (i) of Lemma 3.4,
$\{y_1\}*Y_2\subset M_1$ (note that $y_1\in M_1^\circ$), and thus
$(\{y_1\}*Y_2)\cap N\subset Y_2\subset Y_1*Y_2$ (see footnote 5).

Case 2: $Y_2\subset M_2^\circ$. In this case, by (iii) of Lemma 3.4,
$[y_1y_2]_M\cap N$ contains only one point for any $y_2\in Y_2$.
Then like proving Lemma 2.1, we can prove that $N\cap
(\{y_1\}*[y_2y_2'])$ is a geodesic in
$\{y_1\}*[y_2y_2']\subset\{y_1\}*Y_2$ for any $[y_2y_2']\subset Y_2$
(hint: $\{y_1\}*[y_2y_2']$ can be isometrically embedded into the
unit sphere). This implies that $(\{y_1\}*Y_2)\cap N$ is convex in
$\{y_1\}*Y_2$, so $X:=\{x|x\in[y_1z] \text{ with } z\in
(\{y_1\}*Y_2)\cap N\}$ and $\overline{X^c}$ are  convex in
$S(Y_2)=\{y_1,\bar y_1\}*Y_2$. That is, $S(Y_2)=X\cup_{\partial
X}\overline{X^c}$ with $X, \overline{X^c}\in\mathcal{A}(1)$, which
contradicts Corollary 0.1.

Case 3: $Y_2\cap M_2^\circ\neq\emptyset$ and $Y_2\cap
N\neq\emptyset$. In this case, we select $y_2\in Y_2\cap N$ and
$y_2'\in Y_2\cap M_2^\circ$. Now we consider a geodesic
$[y_2y_2']_{M_2}$. Like proving ``$[y_1p_1]_{M_1}$ has to belong to
$Y_1$'' in the proof of (iii) of Lemma 3.4, we can first conclude
that $[y_2y_2']_{M_2}\subset Y_2$ which implies that
$(\uparrow_{y_2}^{y_2'})_{M_2}\in(\Sigma_{y_2}M_2)^\circ\cap\Sigma_{y_2}Y_2$;
and then by applying the induction on
$\Sigma_{y_2}M=(\Sigma_{y_2}Y_2)*Y_1$ (see 3.2) we have that
$(\uparrow_{y_2}^{y_1'})_M\in\partial(\Sigma_{y_2}M_i)$ for any
$y_1'\in Y_1$. On the other hand, since $y_2\in N$ and $y_1\in
M_1^\circ$, by (i) of Lemma 3.4 $[y_2y_1]_M=[y_2y_1]_{M_1}$. This
implies that $(\uparrow_{y_2}^{y_1})_M\in (\Sigma_{y_2}M_1)^\circ$
(footnote 5), which contradicts
``$(\uparrow_{y_2}^{y_1})_M\in\partial(\Sigma_{y_2}M_i)$''.

\vskip 1mm

Note that in  Cases 2 and 3 we both obtain contradictions, so
$(\{y_1\}*Y_2)\cap N$ has to belong to $Y_2\subset Y_1*Y_2$ (see
Case 1). \hfill $\Box$

\vskip2mm

\noindent{\it Proof of Lemma 3.5}.

We first prove  that $N\cap Y_j\neq\emptyset$. In fact if $N\cap
Y_j=\emptyset$ for $j=1$ or 2, say $j=1$, then we can assume that
$Y_1\subset M_1^\circ$ (or $M_2^\circ$). According to Sublemma
3.5.1, $\{y_1\}*Y_2\subset M_1$ for all $y_1\in Y_1$, i.e.
$M=Y_1*Y_2\subseteq M_1$; a contradiction.

\vskip1mm

As for the latter part of the lemma, we first prove that $Y_j\cap
M_i^\circ\neq\emptyset$ ($i=1, 2$) for $j=1$ or 2. If this is not
true, then either $Y_1, Y_2\subset N$, or $Y_1\subset M_i$ and
$Y_1\cap M_i^\circ\neq\emptyset$ for $i=1$ or 2. Note that in the
latter case, we have $Y_2\subset M_i$ by Sublemma 3.5.1. That is, in
any case we have $Y_1, Y_2\subset M_i$ for $i=1$ or 2. Then by Lemma
3.4, $M=Y_1*Y_2\subseteq M_i$, which is impossible. Now we can
assume that $Y_1\cap M_i^\circ\neq\emptyset$ ($i=1, 2$), and it
remains to show that $Y_2\subset N$. Since $Y_1\cap
M_i^\circ\neq\emptyset$ for $i=1$ and 2, by Sublemma 3.5.1
$Y_2\subset M_1$ and $Y_2\subset M_2$, i.e. $Y_2\subset M_1\cap
M_2=N$. \hfill $\Box$

\vskip2mm

{\noindent \bf 3.6\ \ The proof of Theorem C by assuming $\diam(Y_1*Y_2)<\pi$}

\vskip2mm

According to Lemma 3.5, without loss of generality we assume that
$Y_1\cap M_i^\circ\neq\emptyset$ for $i=1$ and 2, and $Y_2\subset
N$. Let $X_i=Y_1\cap M_i$ for $i=1$ and 2. It follows from (i) and
(ii) of Lemma 3.4 that $$M_i=\bigcup_{x\in X_i}\{x\}*Y_2  \text{ for
} i=1 \text{ and } 2.\eqno{(3.7)}$$ {\bf Claim}: {\it $X_i$ is
convex in $M_i$, and thus $X_i\in\mathcal{A}(1)$}. In order to see
the claim, it suffices to show that there exists $[xx']_{M_i}$ which
belongs to $X_i$ for any $x,x'\in X_i$. We first prove that any
$[xx']_{M_i}\subset X_i$ if $x$ or $x'$, say $x'$, belongs to
$M_i^\circ$. Since $x'\in M_i^\circ$,
$[xx']_{M_i}\setminus\{x\}\subset M_i^\circ$ is a local geodesic in
$M$. If $[xx']_{M_i}\not\subset X_i$ (i.e. $[xx']_{M_i}\not\subset
Y_1$), then by Corollary 3.3.3 (and Remark 3.3.4) $[xx']_{M_i}$ is
of length $\pi$ (so $|xx'|_i=\pi$). Then $M_i=\{x,x'\}*A_i$ for some
$A_i\in \mathcal{A}(1)$ with $\partial M_i=\{x,x'\}*\partial A_i$
(see Corollary A.4.1 in Appendix) which contradicts ``$x'\in
M_i^\circ$''. Next we prove that there exists $[xx']_{M_i}\subset
X_i$ for any $x, x'\in X_i\cap N$. Let $\bar x\in X_i\cap
M_i^\circ$. We have proved that any $[x'\bar x]_{M_i}\subset X_i$.
Moreover, $[x'\bar x]_{M_i}\setminus\{x'\}\subset M_i^\circ$ (see
footnote 5), so $[x\tilde x]_{M_i}\subset X_i$ for any $\tilde
x\in[x'\bar x]^\circ_{M_i}$. Then as $\tilde x$ converges to $x'$,
$[x\tilde x]_{M_i}$ converges to a geodesic $[xx']_{M_i}\subset
X_i$. Hence, the claim is verified.

Note that (3.7) and the above claim imply that $M_i=X_i*Y_2$ and
$\partial M_i=(\partial X_i)*Y_2$. Since $\partial M_1\cong\partial
M_2$ and $M_1\cup_{\partial M_i}M_2=Y_1*Y_2$, $\partial
X_1\cong\partial X_2$ and $Y_1=X_1\cup_{\partial X_i} X_2$. That is,
the proof is done ({\it which together with Lemma 3.1 implies that
the proof of Theorem C is finished}). \hfill $\Box$

%%%%%%%%%%%%%%%%%%%%%%%%%%%%%%%%%%%% Appendix   %%%%%%%%%%%%%%%%%%%%%%%%%%%%%%%%%%%%%%%%

\vskip8mm

\noindent{\large \bf Appendix\ \ \ On joins}

\vskip3mm

{\noindent \bf A.1\ \ On $X*Y$ with $\dim(X)=0$}

\vskip2mm

In $X*Y$, if $\dim(X)=0$, then we make a convention that $X$
consists of either two points with distance equal to $\pi$ or only
one point. In the former case, $X*Y$ is the suspension $S(Y)$; in
the latter case, $X*Y$ is a half suspension.

\vskip2mm

{\noindent \bf Proposition A.1} {\it Let $M\in \mathcal{A}^n(1)$,
and let $p,q\in M$. Then $M=\{p,q\}*X$ for some $X\in
\mathcal{A}^{n-1}(1)$ if and only if $|pq|=\pi$.}

\vskip1mm

\noindent{\it Proof}. The ``only if''  follows from the definition
of the metric of suspensions. We give a brief proof for the ``if''.
{\bf Claim}: {\it Any triangle $\triangle px_1x_2$ with $x_i\neq q$
is isometric to its comparison triangle.} Since $|pq|=\pi$,
$|px_i|+|x_iq|=\pi$ (cf. [BGP]), i.e. any $[px_i]\cup[x_iq]$ is a
geodesic between $p$ and $q$. This implies that $\angle px_1x_2+\angle qx_1x_2=\pi$, and there is a unique geodesic
between $p$ (resp. $q$) and $x_i$. Moreover, ``$|px_i|+|x_iq|=\pi$'' implies
that $\angle \tilde p\tilde x_1\tilde x_2+\angle \tilde q\tilde
x_1\tilde x_2=\pi$, where the angles are in the comparison triangles
of $\triangle px_1x_2$ and $\triangle qx_1x_2$. Note that $\angle px_1x_2\geq \angle \tilde p\tilde x_1\tilde
x_2$ and $\angle qx_1x_2\geq \angle \tilde q\tilde x_1\tilde x_2$ (by the TCT), so
$$\angle px_1x_2=\angle \tilde p\tilde x_1\tilde x_2,\ \angle qx_1x_2=\angle \tilde q\tilde x_1\tilde x_2.$$
Then by the TCT for ``='', $\triangle px_1x_2$ and $\triangle
qx_1x_2$ are isometric to their comparison triangles respectively
(note that there is a unique geodesic between $p$ (resp. $q$) and
$x_2$), i.e. the claim is verified.

Let $X:=\{x\in M||px|=\frac\pi2\}$. By the above claim any triangle
$\triangle px_1x_2$ with $x_i\in X$ is isometric to its comparison
triangle. Hence, $|px|=\frac\pi2$ for any $x\in [x_1x_2]\subset \triangle px_1x_2$, and thus
$[x_1x_2]\subset X$. That is, $X$ is convex in $M$, so $X\in
\mathcal{A}^{n-1}(1)$. Moreover, there is a unique geodesic between $p$
(resp. $q$) and any $x\in X$; and for any $y\in
M$, $y\in [px]$ or $[qx]$, where $x$ is the middle point of the
geodesic $[pq]=[py]\cup[yq]$. Hence, we conclude that $M=\{p,q\}*X$ (see Remark A.3.4 below). \hfill $\Box$

\vskip2mm

{\noindent \bf Remark A.1.1} From the proof of Proposition A.1, we can conclude that {\it if
$|pq|=\pi$ in $M\in \mathcal{A}(1)$, then any triangle $\triangle
pxy\subset M$ with $|px|+|py|+|xy|<2\pi$ is isometric to its
comparison triangle}.

\vskip2mm

{\noindent \bf A.2\ \ An explanation to the definition of the metric
of $X*Y$}

\vskip2mm

In [BGP], the metric of $X*Y$ (see Sec. 0) is given directly. Here
we supply an explanation to its definition, from which we can see
some basic properties of the join.

\vskip2mm

On $X\times Y\times[0,\frac\pi2]/\sim$, where $(x_1,y_1,
a_1)\sim(x_2,y_2, a_2)\Leftrightarrow a_1=a_2=0$ and $x_1=x_2$ or
$a_1=a_2=\frac\pi2$ and $y_1=y_2$, we first use the cosine law of
$\Bbb S^2$ to define
$$\begin{aligned}&\cos|q_1q_2|=\cos a_1\cos a_2+\sin a_1\sin a_2\cos|y_1y_2|,\\
& \cos|r_1r_2|=\cos(\frac\pi2-a_1)\cos(\frac\pi2- a_2)+\sin(\frac\pi2-a_1)\sin(\frac\pi2-a_2)\cos|x_1x_2|,
\end{aligned}\eqno{(\text{A}1)}$$ where $q_i=[(x, y_i, a_i)]$ and $r_i=[(x_i, y, a_i)]$.

Now let $p_i=[(x_i,y_i,a_i)]$ with $i=1$ and 2. Due to (A1),
$\{[(x_i,y_i,t)]|t\in [0,\frac\pi2]\}$ is a geodesic of length
$\frac\pi2$, and for any geodesic $[x_1x_2]\subset X$,
$\{[(x,y_2,t)]|x\in[x_1x_2], t\in [0,\frac\pi2]\}$
($=[x_1x_2]*\{y_2\}$) can be isometrically embedded into $\Bbb S^2$.
Then we can define $|\uparrow_{x_1}^{x_2}\uparrow_{x_1}^{y_2}|$,
$|\uparrow_{x_1}^{x_2}\uparrow_{x_1}^{p_2}|$ and
$|\uparrow_{x_1}^{p_2}\uparrow_{x_1}^{y_2}|$ to be the angles
between the corresponding geodesics in $[x_1x_2]*\{y_2\}\ (\subset
\Bbb S^2)$. It therefore follows that we can embed the four
directions $\uparrow_{x_1}^{x_2}$, $\uparrow_{x_1}^{p_2}$,
$\uparrow_{x_1}^{y_2}$ and $\uparrow_{x_1}^{y_1}$ into $\Bbb S^2$
with
$$|\uparrow_{x_1}^{x_2}\uparrow_{x_1}^{p_2}|+|\uparrow_{x_1}^{p_2}\uparrow_{x_1}^{y_2}|=|\uparrow_{x_1}^{x_2}\uparrow_{x_1}^{y_2}|=\frac\pi2,\ |\uparrow_{x_1}^{x_2}\uparrow_{x_1}^{y_1}|=\frac\pi2,\ |\uparrow_{x_1}^{y_1}\uparrow_{x_1}^{y_2}|=|y_1y_2|.$$
Then we can use the cosine law of $\Bbb S^2$ to define
$|\uparrow_{x_1}^{y_1}\uparrow_{x_1}^{p_2}|$ (i.e. the angle between
$[x_1p_1]$ and $[x_1p_2]$) by
$$\cos|\uparrow_{x_1}^{y_1}\uparrow_{x_1}^{p_2}|=\sin|\uparrow_{x_1}^{x_2}\uparrow_{x_1}^{p_2}|\cos|y_1y_2|,\eqno{(\text{A}2)}$$
and then define $|p_1p_2|$ by
$$\cos|p_1p_2|=\cos a_1\cos|x_1p_2|+\sin a_1\sin|x_1p_2|\cos|\uparrow_{x_1}^{y_1}\uparrow_{x_1}^{p_2}|.\eqno{(\text{A}3)}$$
Because $[x_1x_2]*\{y_2\}$ can be isometrically embedded into $\Bbb S^2$,
$$\begin{aligned}&\cos|x_1p_2|=\cos a_2\cos|x_1x_2|,\\
&\sin a_2=\cos(\frac\pi2-a_2)=\sin|x_1p_2|\cos|\uparrow_{x_1}^{p_2}\uparrow_{x_1}^{y_2}|=\sin|x_1p_2|\sin|\uparrow_{x_1}^{x_2}\uparrow_{x_1}^{p_2}|.
\end{aligned}\eqno{(\text{A}4)}$$
Obviously, plugging (A2) and (A4) into (A3), we obtain that
$$\hskip2.5cm\cos|p_1p_2|=\cos a_1\cos a_2\cos|x_1x_2|+\sin a_1\sin a_2\cos|y_1y_2|.\hskip2cm\Box$$

\vskip2mm

{\noindent \bf A.3\ \ A criterion for the join}

\vskip2mm

{\noindent \bf Proposition A.3} {\it Let $M\in\mathcal{A}(1)$
without boundary, and let $X, Y$ be two convex subsets in $M$ (and
thus $X,Y\in \mathcal{A}(1)$). Then $M=X*Y$ if and only if the
following holds:

\noindent{\rm (i)} $\dim(X)+\dim(Y)+1=\dim(M)$;

\noindent{\rm (ii)} $X$ and $Y$ have empty boundary;

\noindent{\rm (iii)} $|xy|=\frac\pi2$ for any $x\in X$ and $y\in Y$;

\noindent{\rm (iv)} There is a unique geodesic between any $x\in X$
and $y\in Y$.}

\vskip2mm

\noindent{\it Proof}. From the definition of the metric the join, the `only if' is almost obvious. As for the `if', it suffices to show that $X*Y$ can be isometrically embedded into $M$ (this implies that $X*Y$ without boundary is convex in $M$, and thus $X*Y=M$ because $\dim(X*Y)=\dim(M)$ (see footnote 6)).
That is, we need to show that
$$\cos|p_1p_2|=\cos a_1\cos a_2\cos|x_1x_2|+\sin a_1\sin a_2\cos|y_1y_2|,\eqno{(\text{\rm A}6)}$$
for any $x_i\in X$, $y_i\in Y$ and $p_i\in[x_iy_i]$ with $|p_ix_i|=a_i$ ($i=1,2$).
In order to prove this, we consider $\Bbb S^3=S_1^1*S_2^1$ (diam$(S_i^1)=\pi$). Select $\tilde x_i\in S_1^1$, $\tilde y_i\in S_2^1$ and $\tilde p_i\in[\tilde x_i\tilde y_i]$ such that $|\tilde x_1\tilde x_2|=|x_1x_2|$, $|\tilde y_1\tilde y_2|=|y_1y_2|$ and $|\tilde p_i\tilde x_i|=a_i$. Note that $$\cos|\tilde p_1\tilde p_2|=\cos a_1\cos a_2\cos|x_1x_2|+\sin a_1\sin a_2\cos|y_1y_2|.$$

We select a geodesic $[x_1x_2]$ in $X$. Because of (iii) and (iv), $|\uparrow_{x_2}^{y_2}\uparrow_{x_2}^{x_1}|=\frac\pi2$ by the first variation formula ([BGP]), and thus by the TCT for ``='' triangle $\triangle y_2x_1x_2$ (with sides $[y_2x_1], [y_2x_2],$ $[x_1x_2]$) is isometric to its comparison triangle.
This implies that $$|\tilde x_1\tilde p_2|=|x_1p_2|, \text{ and similarly } |\tilde y_1\tilde p_2|=|y_1p_2|; \eqno{(\text{\rm A}7)}$$ i.e., $\triangle \tilde x_1\tilde y_1\tilde p_2$ is the comparison triangle of any
$\triangle x_1 y_1 p_2$. {\bf Claim}: {\it There is a triangle  $\triangle x_1 y_1 p_2$ which is isometric to $\triangle \tilde x_1\tilde y_1\tilde p_2$}. Note that the claim implies (A6), so in the rest we only need to verify the claim.

Since the triangle $\triangle y_2x_1x_2$ is isometric to its comparison triangle, $\triangle y_2x_1x_2$ bounds a domain which can be isometrically embedded into $\Bbb S^2$ ([GM]). Select the geodesic $[x_1p_2]$ in this domain. For any $p_3\in [x_1p_2]^\circ$ (resp. $\tilde p_3\in [\tilde x_1\tilde p_2]$ with $|\tilde p_3\tilde x_1|=|p_3x_1|$), there is $x_3\in [x_1x_2]$ (resp. $\tilde x_3\in [\tilde x_1\tilde x_2]$ with $|\tilde x_3\tilde x_1|=|x_3x_1|$)  such that $p_3\in [x_3y_2]$ (resp. $\tilde p_3\in [\tilde x_3\tilde y_2]$ with $|\tilde x_3\tilde p_3|=|x_3p_3|$). Then we can conclude that
$|\tilde y_1\tilde p_3|=|y_1p_3|$ as same as getting (A7). Hence, by the TCT for ``='' there is a triangle  $\triangle x_1 y_1 p_2$ which is isometric to
its comparison triangle, i.e. the above claim holds. \hfill $\Box$

\vskip2mm

From the proof of Proposition A.3 (together with A.2), we can conclude:

\vskip2mm

{\noindent \bf Corollary A.3.1} {\it Let $X*Y$ be the join of $X,Y\in \mathcal{A}(1)$. For any $p_i=[(x_i,y_i,a_i)]\in X*Y$ and any geodesic $[x_1p_2]$, there is a geodesic
$[p_1p_2]$ such that the triangle $\triangle p_1x_1p_2$ composed by $[x_1p_1]$, $[x_1p_2]$ and $[p_1p_2]$ is isometric to its comparison triangle.}

\vskip2mm

{\noindent\bf Remark A.3.2} In Proposition A.3, (iii) and (iv) can be replaced by ``$|xy|=\frac\pi2$ and $[xy]$ is unique for any $x\in X'$ and $y\in Y'$, where $X'$ and $Y'$ are dense in $X$ and $Y$ respectively''. Obviously, due to the continuity of the distance function, $|xy|=\frac\pi2$ for any $x\in X$ and $y\in Y$.
Then similarly, we can prove that $X*Y$ can be isometrically embedded into $M$ (so $X*Y=M$) once we show that
$$\lim_{j\to+\infty} [x_{1j}y_{1j}]=\lim_{j\to+\infty}[x_{2j}y_{2j}]$$
for any $x_{1j},x_{2j}(\in X')\to x$ and $y_{1j}, y_{2j}(\in Y')\to y$ as $j\to+\infty$ with $x\in X$ and $y\in Y$.
By considering geodesics $[x_{1j}y_{2j}]$ in addition, like getting (A7) we conclude that 
$$\angle y_{1j}x_{1j}y_{2j}=|y_{1j}y_{2j}| \text{ and } \angle x_{1j}y_{2j}x_{2j}=|x_{1j}x_{2j}|.$$
It then follows (p.5 of [BGP]) that
$$\hskip3.5cm\lim_{j\to+\infty}[x_{1j}y_{1j}]=\lim_{j\to+\infty}[x_{1j}y_{2j}]=\lim_{j\to+\infty}[x_{2j}y_{2j}].\hskip3cm\Box$$

{\noindent \bf Remark A.3.3} From the proof of Proposition A.3 and Remark A.3.2, we can conclude that: {\it Let $X$ and $Y$ be two convex subsets in $M\in\mathcal{A}(1)$. If $|xy|=\frac\pi2$ and $[xy]$ is unique for any $x\in X'$ and $y\in Y'$, where $X'$ and $Y'$ are dense in $X$ and $Y$ respectively, then $X*Y$ can be isometrically embedded into $M$.} \hfill $\Box$

\vskip2mm

{\noindent \bf Remark A.3.4} From Remark A.3.3, we can conclude that: {\it Let $X$ and $Y$ be two convex subsets in $M\in\mathcal{A}(1)$. If $|xy|=\frac\pi2$ and $[xy]$ is unique for any $x\in X$ and $y\in Y$, and if for any $p\in M$ there exist $x\in X$ and $y\in Y$ such that $p\in[xy]$, then $M=X*Y$.} \hfill $\Box$

\vskip2mm

{\noindent \bf A.4\ \ Direction spaces at the points on the `bottoms' of $X*Y$}

\vskip2mm

{\noindent \bf Proposition A.4} {\it  Let $M=X*Y$ with $X,Y\in \mathcal{A}(1)$. Then for any $x_0\in
X\subset M$ (resp. $Y\subset M$), the direction space
$\Sigma_{x_0}M=\Sigma_{x_0}X*Y$ (resp. $X*\Sigma_{x_0}Y$). }

\vskip2mm

Recall that $p\in M\in \mathcal{A}(k)$ is a boundary point if $\Sigma_pM (\in \mathcal{A}(1))$ has nonempty
boundary ([BGP]). Hence, Proposition A.4 has an immediate corollary.

\vskip2mm

{\noindent \bf Corollary A.4.1} {\it  Let $M=X*Y$ with $X,Y\in \mathcal{A}(1)$. Then $M$ has nonempty boundary if and only if at least one of $X$ and $Y$ has nonempty boundary,
and $\partial M=(\partial X*Y)\cup (X*\partial Y)$. }

\vskip2mm

\noindent{\it Proof of Proposition A.4}. Since $|xy|=\frac\pi2$ and $[xy]$ is unique
for any $x\in X$ and $y\in Y$ (see Proposition A.3), by the TCT for
``='' the map $$i:Y\to \Sigma_{x_0}M \text{ defined  by } y\mapsto
\uparrow_{x_0}^y$$ is an isometrical embedding; and by the first
variation formula ([BGP]) $|\uparrow_{x_0}^y\xi|=\frac\pi2$ (in
$\Sigma_{x_0}M$) for any $y\in Y$ and $\xi\in\Sigma_{x_0}X$. 
Note that $i(Y)$ ($\stackrel{\rm{iso}}{\cong} Y$) is convex in $\Sigma_{x_0}M$ because $i$ is an isometrical
embedding; and $\Sigma_{x_0}X$ is convex in $\Sigma_{x_0}M$ because $X$ is convex in $M$. {\bf
Claim 1}: {\it there is a unique geodesic between $\uparrow_{x_0}^y$
and any $\xi\in (\Sigma_{x_0}X)'$.} Recall that
$(\Sigma_{x_0}X)'=\{\xi\in\Sigma_{x_0}X| \exists\ [x_0x]\subset X
\text{ s.t. }\xi=\uparrow_{x_0}^{x}\}$, which is dense in
$\Sigma_{x_0}X$ ([BGP]). Then by Remark A.3.3, $\Sigma_{x_0}X*Y$ can
be isometrically embedded into $\Sigma_{x_0}M$; and by Corollary A.3.1, it is not hard to see that
$(\Sigma_{x_0}M)'\subseteq(\Sigma_{x_0}X)'*Y$. It then follows that $\Sigma_{x_0}M=\Sigma_{x_0}X*Y$.

In the rest,  we only need to verify Claim 1. Select $[x_0x]\subset
X$ such that $\xi=\uparrow_{x_0}^{x}$. Since $[x_0x]*\{y\}\subset M$
can be isometrically embedded into $\Bbb S^2$ (see (A.2)), $\gamma:=\{\uparrow_{x_0}^z|z\in [xy] \text{ and
} [x_0z]\subset [x_0x]*\{y\}\}$ is a
geodesic between $\uparrow_{x_0}^{y}$ and $\xi$ in $\Sigma_{x_0}M$.

{\bf Claim 2}: {\it For any $[x_0p]\subset M$,
$|\uparrow_{x_0}^{x}\uparrow_{x_0}^p|+|\uparrow_{x_0}^p\uparrow_{x_0}^{y}|\geq\frac{\pi}{2}$,
and the ``$=$'' holds only if $\uparrow_{x_0}^p\in \gamma$.} Let $p=[(x',y',t)]\in X*Y$, and select a geodesic
$[x_0x']\subset X$ such that $[x_0p]\subset [x_0x']*\{y'\}$. Note that we have proved that
$|\uparrow_{x_0}^{x'}\uparrow_{x_0}^{y}|=|\uparrow_{x_0}^{x'}\uparrow_{x_0}^{y'}|=|\uparrow_{x_0}^{x}\uparrow_{x_0}^{y'}|=\frac{\pi}{2}$, and 
$$|\uparrow_{x_0}^{x'}\uparrow_{x_0}^p|+|\uparrow_{x_0}^p\uparrow_{x_0}^{y'}|=|\uparrow_{x_0}^{x'}\uparrow_{x_0}^{y'}|=\frac{\pi}{2},$$
i.e. $\uparrow_{x_0}^p$ belongs to a geodesic $[\uparrow_{x_0}^{x'}\uparrow_{x_0}^{y'}]\subset\Sigma_{x_0}M$. By considering triangles
$\triangle\uparrow_{x_0}^{y'}\uparrow_{x_0}^{x'}\uparrow_{x_0}^{x}$ and $\triangle\uparrow_{x_0}^{x'}\uparrow_{x_0}^{y'}\uparrow_{x_0}^{y}$ (both of which contain the side $[\uparrow_{x_0}^{x'}\uparrow_{x_0}^{y'}]$) and their
comparison triangles, it is not hard to see that
$$|\uparrow_{x_0}^p\uparrow_{x_0}^{x}|\geq|\uparrow_{x_0}^p\uparrow_{x_0}^{x'}|\text{ and } |\uparrow_{x_0}^p\uparrow_{x_0}^{y}|\geq|\uparrow_{x_0}^p\uparrow_{x_0}^{y'}|;\eqno{(\text{A}8)}$$
moreover, ``the two ``$=$'' hold'' implies that  $|yy'|=|\uparrow_{x_0}^{x}\uparrow_{x_0}^{x'}|=0$,
i.e. $y=y'$ and $[x_0x]\subseteq[x_0x']$ or vice versa, which implies that  $\uparrow_{x_0}^p\in \gamma$ (note that $[x_0p]\subset [x_0x']*\{y'\}$).
Hence, Claim 2 follows.

Due to Claim 2, if Claim 1 is not true, i.e. there is another
geodesic $[\uparrow_{x_0}^{y}\uparrow_{x_0}^{x}]\neq\gamma$ in
$\Sigma_{x_0}M$, then
$[\uparrow_{x_0}^{y_0}\uparrow_{x_0}^{x_1}]^\circ\cap(\Sigma_{x_0}M)'=\emptyset$.
Select the middle points of $\gamma$ and
$[\uparrow_{x_0}^{y}\uparrow_{x_0}^{x}]$, denoted by
$\uparrow_{x_0}^{p}$ and $\eta$ respectively. For convenience, we
denote $|\uparrow_{x_0}^{p}\eta|$ to be $C$. Recall that for any $0<\delta<<1$, there exists
$[x_0p_1]\subset M$ with $p_1=[(x_1,y_1,t)]$ such that
$|\uparrow_{x_0}^{p_1}\eta|<C\delta$ ([BGP]), and thus
$$\left||\uparrow_{x_0}^{p_1}\uparrow_{x_0}^{y}|-\frac\pi4\right|<C\delta
\text{ and }
\left||\uparrow_{x_0}^{p_1}\uparrow_{x_0}^{x}|-\frac\pi4\right|<C\delta.\eqno{(\text{A}9)}$$
Then similar to getting (A8), by analyzing triangles
$\triangle\uparrow_{x_0}^{y_1}\uparrow_{x_0}^{x_1}\uparrow_{x_0}^{x}$
and
$\triangle\uparrow_{x_0}^{x_1}\uparrow_{x_0}^{y_1}\uparrow_{x_0}^{y}$
(both of which contain the side $[\uparrow_{x_0}^{x_1}\uparrow_{x_0}^{y_1}]\ni
\uparrow_{x_0}^{p_1}$) and their comparison triangles, we can
conclude that
$$|\uparrow_{x_0}^{x}\uparrow_{x_0}^{x_1}|<C_1\sqrt\delta  \text{ and } |yy_1|<C_1\sqrt\delta \text{\ \ for some constant  } C_1. \eqno{(\text{A}10)}$$
Now we consider $[x_0x]*\{y_1\}$, and select $[x_0p_2]\subset [x_0x]*\{y_1\}$ with $p_2\in[xy_1]$ and
$|\uparrow_{x_0}^{y_1}\uparrow_{x_0}^{p_2}|=|\uparrow_{x_0}^{p_2}\uparrow_{x_0}^{x}|$
$=\frac\pi4$. Since  $y_1\stackrel{\delta\to0}{\longrightarrow}y$
(see (A10)), $[x_0x]*\{y_1\}\longrightarrow[x_0x]*\{y\}$ as $\delta\to0$. Hence,
$|\uparrow_{x_0}^{p}\uparrow_{x_0}^{p_2}|<\chi(\delta)$, where
$\chi(\delta)\stackrel{\delta\to0}{\longrightarrow}0$, and thus
$$\left||\uparrow_{x_0}^{p_2}\uparrow_{x_0}^{x}|-\frac{\pi}{4}\right|<\chi(\delta)
\text{ and
}\left||\uparrow_{x_0}^{p_2}\uparrow_{x_0}^{p_1}|-C\right|<\chi(\delta).\eqno{(\text{A}11)}$$
Without loss of generality, we can assume that
$|x_0x|=|x_0x_1|=\epsilon$. Due to (A9) and (A11), in the triangles
$\triangle p_2xx_0$ and $\triangle p_1x_1x_0$,
$$\left||xp_2|-\epsilon\right|<\chi(\delta,\epsilon)\epsilon
\text{ and }
\left||x_1p_1|-\epsilon\right|<\chi(\delta,\epsilon)\epsilon,\eqno{(\text{A}12)}$$
and
$$\left||x_0p_2|-\sqrt2\epsilon\right|<\chi(\delta,\epsilon)\epsilon
\text{ and }
\left||x_0p_1|-\sqrt2\epsilon\right|<\chi(\delta,\epsilon)\epsilon,$$
where $\chi(\delta,\epsilon)\to0$ as $\delta\to 0$ and $\epsilon\to
0$. Then due to (A10) and (A11), in the triangles $\triangle
x_0xx_1$ and $\triangle x_0p_1p_2$
$$|xx_1|<(1+\chi(\delta,\epsilon))C_1\sqrt\delta\epsilon
\text{ and } \left||p_1p_2|-2\sqrt2\sin\frac
C2\epsilon\right|<\chi(\delta,\epsilon)\epsilon.\eqno{(\text{A}13)}$$
On the other hand, since any triangle $\triangle y_1xx_1$ with
$|y_1x|=|y_1x_1|=\frac\pi2$ is isometric to its comparison triangle
(Corollary A.3.1), (A12) implies that
$$\left|\frac{|p_1p_2|}{|xx_1|}-1\right|<C_2\epsilon^2\text{\ \
for some constant  } C_2,$$ which contradicts (A13). Note that this
contradiction is drawn under the assumption that Claim 1 is not true.
Hence, Claim 1 holds, and the proof is done. \hfill $\Box$

\vskip2mm

{\noindent \bf Remark A.4.2}  Once we finished the proof of Proposition A.4,
we can say that any $\triangle\uparrow_{x_0}^{y'}\uparrow_{x_0}^{x'}\uparrow_{x_0}^{x}$ and $\triangle\uparrow_{x_0}^{x'}\uparrow_{x_0}^{y'}\uparrow_{x_0}^{y}$
in the proof of Proposition A.4 (before (A8)) are isometric to their
comparison triangles respectively. Consequently, (A8) has a precise formulation
$$\begin{aligned}&\cos|\uparrow_{x_0}^{p}\uparrow_{x_0}^{x}|=\cos|\uparrow_{x_0}^{p}\uparrow_{x_0}^{x'}|\cos|\uparrow_{x_0}^{x}\uparrow_{x_0}^{x'}|,\\
&\cos|\uparrow_{x_0}^{p}\uparrow_{x_0}^{y}|=\cos|\uparrow_{x_0}^{p}\uparrow_{x_0}^{y'}|\cos|yy'|.
\end{aligned}$$

\vskip2mm

{\noindent \bf A.5\ \ On the diameter of the join}

\vskip2mm

{\noindent \bf Proposition A.5}
{\it Let $X*Y$ be the join of $X,Y\in \mathcal{A}(1)$.
Then $$\diam(X*Y)=\max\{\diam(X), \diam(Y), \frac\pi2\}\footnote{This part of the proposition was mentioned in [GP].}.$$
Moreover, if $\diam(X*Y)>\frac\pi2$ and $\diam(X)>\diam(Y)$, and if $|pq|=\diam(X*Y)$, then $p,q\in X$.}

\vskip2mm

Proposition A.5 is an immediate corollary of the following lemma.

\vskip2mm

{\noindent \bf Lemma A.5.1}
{\it Let $X*Y$ be the join of $X,Y\in \mathcal{A}(1)$, and let $p_i=[(x_i, y_i, t_i)]\in X*Y$ ($i=1,2$).
If $|p_1p_2|=\diam(X*Y)>\frac\pi2$, then either $|x_1x_2|=|y_1y_2|=|p_1p_2|$ and $t_1=t_2$, or $|x_1x_2|=|p_1p_2|>|y_1y_2|$ and $t_1=t_2=0$, or $|y_1y_2|=|p_1p_2|>|x_1x_2|$ and $t_1=t_2=\frac\pi2$.}

\vskip1mm

\noindent{\it Proof}. We first assume that $|x_1x_2|\leq |y_1y_2|$. Then by the definition of the metric of joins,
$$\begin{aligned} \cos|p_1p_2|=&\cos t_1\cos t_2\cos|x_1x_2|+\sin t_1\sin t_2\cos|y_1y_2|\\
=& \cos t_1\cos t_2(\cos|x_1x_2|-\cos|y_1y_2|)+\cos(t_1-t_2)\cos|y_1y_2|.
\end{aligned}$$
Since $|p_1p_2|>\frac\pi2$ and $|p_1p_2|\geq|y_1y_2|$, it is clear that
$$\cos t_1\cos t_2(\cos|x_1x_2|-\cos|y_1y_2|)=0,\ t_1-t_2=0\text{ and } |y_1y_2|=|p_1p_2|$$
(note that $\cos t_1\cos t_2(\cos|x_1x_2|-\cos|y_1y_2|)\geq0$ and $\cos(t_1-t_2)\geq0$).
That is either $|x_1x_2|=|y_1y_2|=|p_1p_2|$ and $t_1=t_2$, or $|y_1y_2|=|p_1p_2|>|x_1x_2|$ and $t_1=t_2=\frac\pi2$. Similarly, if $|x_1x_2|>|y_1y_2|$, then $|x_1x_2|=|p_1p_2|>|y_1y_2|$ and $t_1=t_2=0$.
\hfill $\Box$

%%%%%%%%%%%%%%%%%%%%%%%%%%%%%%%%%%%%   reference   %%%%%%%%%%%%%%%%%%%%%%%%%%%%%%%%%%%%%%%

\noindent School of Mathematical Sciences (and Lab. math. Com.
Sys.), Beijing Normal University, Beijing, 100875
P.R.C.\\
e-mail: suxiaole$@$bnu.edu.cn; wwyusheng$@$gmail.com

\vskip2mm

\noindent Mathematics Department, Capital Normal University,
Beijing, 100037 P.R.C.\\
e-mail: hwsun$@$bnu.edu.cn

\end{document}